\documentclass[a4 paper,12pt]{article}
\usepackage[T1]{fontenc}
\usepackage[latin1]{inputenc}
\usepackage[german,english]{babel}
\usepackage{graphicx}
\usepackage{mathrsfs}
\usepackage{amsfonts}
\usepackage{pifont}
\usepackage{amsxtra}
\usepackage{amssymb}\usepackage{extarrows} 
\usepackage{eufrak}
\usepackage{pst-all}
\usepackage{multido}
\usepackage{amsmath}
\usepackage{color}
\usepackage{makeidx}
\usepackage{multicol}
\usepackage[absolut]{overpic}
\usepackage{hyperref} 
\usepackage[rm,bf,tiny,center]{titlesec}
\titlelabel{\thetitle.\enspace}
\newcommand{\ger}[1]{\mathfrak{#1}}

\newcommand{\sma}{\triangleleft}

\newcommand{\OO}{\mathcal{O}}

\newcommand{\ma}{\leqslant}

\newcommand{\CC}{\mathbb{C}}

\newcommand{\NN}{\mathbb{N}}

\DeclareMathOperator{\id}{id}

\DeclareMathOperator{\Inv}{Inv}
\DeclareMathOperator{\Red}{R_{\textbf{d}}}
\DeclareMathOperator{\Bruh}{\leqslant}
\newtheorem{numeris}{}[section]
\newtheorem{subnum}{}[numeris]

\usepackage{geometry}

\title{Permutations  with Verma Multiplicities $[M(p):L(q)] \geq 2$}
\author{Daiva Pu\v{c}inskait\.{e}}
\date{\today}

\usepackage{fancyhdr}
\begin{document}
\begin{center}
\textsc{Permutations  with Verma Multiplicities} $[M(p):L(q)] \geq 2$
\end{center}
\begin{center}
Daiva Pu\v{c}inskait\.{e}
\end{center}

\begin{small}
\begin{center}
\begin{minipage}{13.4cm}
\textsc{Abstract.}  We consider permutations $q$ in the symmetric group, $S_n$,  whose Verma multiplicities in the principal block of $\mathcal{O}(\mathfrak{sl}_n)$ satisfy  $[M(p):L(q)] \geq 2$. We present a construction along with a diagrammatic visualization,  showing how  permutations in $S_{n}$ with this property, generates a family of permutations in $S_{n+1}$ that share the same multiplicity property.  While the method does not recover all such permutations in $S_{n+1}$, it systematically generates many new examples. In addition, we present all permutations in $S_{5}$ with $[M(\id_n):L(q)] \geq 2$ in a Bruhat diagram and and describe describe all  permutations in  $S_6$ and $S_7$ with non-simple Verma multiplicities.
\end{minipage}
\end{center}\end{small}

\begin{center}
\textbf{Introduction}
\end{center}

The Bernstein-Gelfand-Gelfand category  $\OO(\ger{g})$ of a simple, finite dimensional, complex Lie Algebra $\ger{g}$ decomposes as a  direct sum of indecomposable subcategories (called blocks)
$\OO_{\lambda}(\ger{g})$, parametrized by   certain weights $\lambda$ (see \cite{BGG}). The  essential modules in the principle block  $\OO_{0}(\ger{g})$    (like  Verma modules, simple modules, projective indecomposable modules, etc.) are  parameterized  by the elements of  the Weyl group $W$ of $\mathfrak{g}$.  
Particular attention in the 
area of the representation theory of $\OO_{0}(\ger{g})$ is paid to the  composition factor multiplicities $[M(w):L(u)]$  of a simple module $L(u)$ in  a  Verma 
 module $M(w)$, which  is strongly related to the Bruhat order, $(W, \ma)$.  The block $\OO_{0}(\mathfrak{g})$ is equivalent to the module category of a finite-dimensional quasi-hereditary algebra $A_{0}(\mathfrak{g})$, whose structure fundamentally depends on multiplicities $[M(\mathrm{id}):L(u)]$ of the principal Verma module $M(\mathrm{id})$ (see \cite{Sorgel}, \cite{CPS}, \cite{Str}).

It is well know
$[M(w),L(u)]=0$ 
if and only if  $w\not< u$  (see \cite{BGG2}).  A criterion developed by Jantzen for $[M(w),L(u)]=1 $  in \cite{Jantzen}   can be reformulated into the following \\

 \textbf{Theorem.} (Multiplicity-One Criterion via Bruhat Intervals) \textit{Let $w , u \in W$  and
$\Lambda^{(w,u)}:=\left\{v\in W\mid w\ma v\ma u\right\}$.  The following
statements are equivalent}
\begin{itemize}
        \item[(i)]  $\left[M(w):L(u)\right]= 1$.
        \item[(ii)]    
				$\left|\left\{v\in \Lambda^{(w,u)}\mid l(v)=l(w)+i
				\right\}\right| = \left|\left\{v\in \Lambda^{(w,u)}\mid l(v)=l(u)-i
\right\}\right|$, for all $i\geq 0$.

\end{itemize}

For more details and information see \cite[8.7]{Humphreys1}.\\[1mm]
\parbox{11.4cm}{\begin{small}
This theorem is visualized in the figure to the right: It is well known
that $\id \in W$  is the uniquely-determined minimal element
of the length  $0$,  and $\omega\in W$  (the longest element) is the 
uniquely-determined maximal element in 
$(W,\ma)$.    The gray  area
represents the subset $\Lambda^{(w,u)}$ of $W$.
We have  $[M(w):L(u)]=1$ if and only if the number of elements in
$\Lambda^{(w,u)}$ of  length $l(w)+i$
 is equal to the number of elements in $\Lambda^{(w,u)}$ of 
length $l(u)-i$.
For example,  for $\ger{g} =\ger{sl}_4(\CC)$ we have $W\cong S_4$. For
the permutation
$\sigma =(4231)$ there are four elements $v\in W$ with $\sigma < v$ and
$l(v)=l(\sigma)-1$, but
 thee  elements $v$ in $W$ with $l(v)=1$. Thus $[M(\id):L(\sigma)]\geq 2$ (see  ~\ref{crit_in_S_4}).
\end{small}
 }  \\[-4.7cm]

   \hspace*{6.7cm}\psset{xunit=0.83cm,yunit=0.83cm,runit=1cm}
\begin{pspicture}(0,0)(0,0)
\rput(8.05,-1.95){
\psclip{
\pscircle[linestyle=none](-0.7,0){1.41cm}
\pscircle[linestyle=none](0.7,0){1.41cm}
        }
        \psframe*[linecolor=gray!30](-2.5,-2.5)(7,7)
\endpsclip
}
\end{pspicture}

\hspace*{6.7cm}\psset{xunit=0.68mm,yunit=0.68mm,runit=1mm}
\begin{pspicture}(0,30)(0,0)
\rput(100,17){
\begin{tiny}
\rput[l](28,19.6){\rnode{oben}{}}
\rput[l](28,0.2){\rnode{obenx1}{}}
\rput[l](28,-41.6){\rnode{unten}{}}
\rput[l](28,-36.4){\rnode{untenx1}{}}
\rput[l](26.5,18.8){\rnode{oben1id}{$- l(\omega)$}}
\rput[l](26.5,9){\rnode{oben1}{$- l(u)$}}
\rput[l](26.5,1){\rnode{oben2}{$- $}}
\rput[l](26.5,-22){\rnode{oben3}{$- $}}
\rput[l](28,-21){\rnode{oben3x}{}}
\rput[l](26.5,-6){\rnode{oben7}{$- \ l(u)-i$}}
\rput[l](26.5,-30){\rnode{oben4}{$- \ l(w)$}}
\rput[l](28,-31){\rnode{oben4x}{}}
\rput[l](26.5,-15){\rnode{oben5}{$- \ l(w)+i $}}
\rput[l](26.5,-41){\rnode{oben6}{$- \ 0 $}}
\rput(0,19.6){\rnode{Aid}{$\omega$}}
\rput(0,9.6){\rnode{A}{$u$}}
\rput(0,-41.6){\rnode{B}{$\id$}}

\rput(20,-17){\rnode{aa}{$\textcolor{white}{a} $}}
\rput(-20,-17){\rnode{bb}{$\textcolor{white}{a} $}}
\rput(-11,14){\rnode{A1}{$$}}
\rput(-7.9,2){\rnode{A11}{$$}}
\rput(7.9,2){\rnode{A21}{$$}}
\rput(-11,-35){\rnode{B1}{}}
\rput(11,14){\rnode{A2}{$$}}
\rput(11,-35){\rnode{B2}{}}
              \rput(0,-30){\rnode{3}{$ w $}}
              \rput(7.9,-23){\rnode{10}{$$}}
              \rput(0,-22){$\cdots$}
               \rput(0,1){$\cdots$}
               \rput(0,-5.6){$\cdots$}
                \rput(0,-15){$\cdots$}
              \rput(-7.9,-23){\rnode{11}{}}

                      \psset{nodesep=1.3pt}
               \ncline[linestyle=dashed]{-}{oben}{unten}   
                 \ncline{-}{oben}{obenx1}   
                 \ncline{-}{unten}{untenx1}  
                 \ncline{-}{oben3x}{oben4x}        
              \ncline{<->}{3}{10}
             \ncline{<->}{3}{11}            
              \ncline{<->}{Aid}{A1}
               \ncline{<->}{A}{A11}
                \ncline{<->}{A}{A21}
              \ncline{<->}{B1}{B}
              \ncline{<->}{Aid}{A2}
              \ncline{<->}{B2}{B}
              \ncarc[linewidth=0.5pt,arcangle=-15,linestyle=dashed]{<->}{sn}{20}
              \ncarc[linewidth=0.5pt,arcangle=15,linestyle=dashed]{<->}{ln}{10}
              \ncarc[arcangle=34,linestyle=dashed]{<->}{A21}{10}
              \ncarc[arcangle=-34,linestyle=dashed]{<->}{A11}{11}
              
              \ncarc[arcangle=-25,linestyle=dashed]{<->}{B}{20}
              \ncarc[arcangle=25,linestyle=dashed]{<->}{B}{21}
               \psset{nodesep=0.3pt}
               
              \ncarc[arcangle=20,linestyle=dashed]{<-}{A2}{aa}
              \ncarc[arcangle=20,linestyle=dashed]{->}{aa}{B2}
              \ncarc[arcangle=-20,linestyle=dashed]{<-}{A1}{bb}
              \ncarc[arcangle=20,linestyle=dashed]{<-}{B1}{bb}
              \rput(0,-15){\rnode{dez-i}{\setlength{\fboxsep}{1.7pt}{\fcolorbox{black}{gray!30}{$\textbf{w}_1 \  \  \  \       \cdots \    \  \   \  \  \textbf{w}_r  $}}}}
             \rput(0,-5.6){\rnode{dez-m-i}{\setlength{\fboxsep}{1.7pt}{\fcolorbox{black}{gray!30}{$\textbf{u}_1 \  \  \  \ \      \cdots \    \  \   \  \  \textbf{u}_r  $}}}}
              \end{tiny}
              }
\end{pspicture} 
\psset{xunit=0.68mm,yunit=0.68mm,runit=1mm}
\begin{pspicture}(0,0)(0,0)
\rput(97,17){
\rput(5,0){\rnode{w1}{}}
 \rput(19,24){\rnode{v1}{$\Lambda^{(w,u)}$}}
              \rput(15,19){\rnode{v2}{}}
 \ncarc[linewidth=1pt,arcangle=25,linecolor=black]{->}{v2}{w1}
}
\end{pspicture}
 \text{ }\\[1.6cm]

 	Kazhdan and Lusztig
described these multiplicities Verma modules in terms of the values at 1 of
certain polynomials related to the Iwahori-Hecke algebra of  the Weyl
group $W$ (see \cite{Kazhdan},\cite{Brenti}). The multiplicity-one criterion can be deduced from the Kazhdan-Lusztig conjecture: since the multiplicity $[M(w):L(u)] = P_{u,w}(1)$, the condition $[M(w):L(u)] = 1$ is equivalent to $P_{u,w}(q) = 1$, which follows from the positivity and symmetry properties of Kazhdan-Lusztig polynomials in finite Weyl groups. 

The study of Verma multiplicities is fundamental in category $\OO$ and provides an important link between the representation theory of Lie algebras, Kazhdan-Lusztig theory, and the combinatorics of the Bruhat order.

In   $\OO_{0}(\ger{g})$ we have $[M(w):L(u)] \geq  2$ if and only if $w<u$ and $[M(w):L(u)] \neq  1$. Combining this with "Multiplicity One Criteria" we have that 
$[M(w):L(u)] \geq  2$ if and only if  
				$\left|\left\{v\in \Lambda^{(w,u)}\mid l(v)=l(w)+i
				\right\}\right| \neq \left|\left\{v\in \Lambda^{(w,u)}\mid l(v)=l(u)-i
\right\}\right|$, for some $i\geq 0$.

In this paper, we are focusing on the principle block of the BGG category $\OO(\ger{sl}_{n}(\CC))$. Since the Weyl group of $\ger{sl}_{n}(\CC)$ is isomorphic to  the symmetric group $S_{n}$, and  for all $u,w\in S_n$ with $w\ma u$ (with respect to Bruhat order) we have   $M(u)\subseteq M(w)$, so   $[M(u),L(q)]\geq  2$ implies $[M(w),L(q)]\geq  2$. Thus the study of composition factors of Verma modules we will reduces to analyzing \textit{principle} Verma module $M(\id_{n})$, since for  the identity element $\id_{n}\in S_{n}$ we have $\id_n\ma w$ for all $w\in S_n$. Permutations $q$ with $[M(\id_{n}),L(q)]\geq  2$ indicate  increased complexity in the structure  of Verma modules, marking cases when local submodules of $M(\id_{n})$ are not Verma modules. Accordingly, we refer to such permutations as \textit{critical permutation} throughout this paper.

 Since $S_{n+1}$ is the disjoint union of its left (or right) cosets of $S_n$, we study how permutations $q \in S_n$ with $[M(\mathrm{id}_n):L(q)] \geq 2$ induce permutations in $S_{n+1}$ lying in the associated cosets, and how these inherit the multiplicity property $[M(\mathrm{id}_{n+1}):L(q')] \geq 2$. For a fixed $n\in \NN$ let
\begin{center}
$
\mathcal{C}_{n} := \left\{ q \in S_{n} \,\middle|\, [M(\id_{n}):L(q)] \geq 2 \right\},
$
\end{center}
  $c_{i}:= s_{i}\cdot s_{i+1} \cdots s_{n}\in S_{n+1}$, and  $d_{i}:=\left(c_{i}\right)^{-1}$, where $s_j$ denotes the adjacent transposition swapping $j$ and $j+1$. 
The main result of this paper is the Theorem~\ref{Main_Theorem}, whose key assertion is stated as follows. \\

\textbf{Theorem} (Lifting Multiplicity $\geq 2$). \textit{If} $q\in \mathcal{C}_{n} $, \textit{then} $c_{i}\cdot q $ \textit{and} $  q\cdot d_{i}\in \mathcal{C}_{n+1} $,   \textit{for each} $i$. \\

Since $ [M(\id_n):L(p)] = [M(\id_n):L(p^{-1})] = [M(\id_n):L(\omega_n \cdot p \cdot \omega_n)]$ where $\omega_n$ is the longest element of $S_{n}$,  the Theorem above yields numerous permutations in  $ \mathcal{C}_{n+1} $  from those in  $\mathcal{C}_{n} $. 
In this paper, we present a detailed proof for permutations  satisfying 
\begin{center}
$\left|\left\{p\in S_{n}\mid  p\ma q \ and  \ l(p)=1 
				\right\}\right| \neq \left|\left\{p\in S_{n}\mid p\ma q \ and  \ l(p)=l(q)-1\right\}\right|  $
\end{center}
solely through the structural properties of permutations and their relations with respect to the Bruhat order,  leaving the general case for future work.

It should be noted that not all permutations in $\mathcal{C}_{n+1}$ arise from elements of $\mathcal{C}_{n}$   via the construction described in the theorem. It is well known that $ |\mathcal{C}_4| = 2 $ and $ |\mathcal{C}_5| = 32 $ (see \cite[5.24]{Jantzen}). Notably, 26 of the permutations in $ \mathcal{C}_5 $ arise from those in $ \mathcal{C}_4 $ via the construction described in the Theorem. 
In this paper, we also identify    permutations that  in $\mathcal{C}_{n+1}$  that do not originate from $\mathcal{C}_{n}$ in this way, yet possess a specific structure and satisfy $[M(\id_{n+1}):L(q)] \geq 2 $.

The paper is organized as follows. Sections 1 and 2 provide brief descriptions of well-known concepts related to symmetric groups that will be used throughout the paper, and 
 notations, definitions related to \textit{critical permutations} as well as a brief account of  elements  in $\mathcal{C}_{n}$ for $n\leq 7$. In Section 3, we present a graphical visualization that facilitate understanding of the Theorem~\ref{Main_Theorem}. Furthermore, we fully characterize the permutations in $\mathcal{C}_{5}$ in 
locate them in the symbolic Bruhat diagram of $S_5$. For $n=6,7$ we describe all permutations that are not `\textit{inherited}' from $\mathcal{C}_{n-1}$, as determined by the method of Theorem~\ref{Main_Theorem}. Section 4 analyzes the structure of permutations, focusing on reduced decompositions and free-fall inversions that are  essential for the proof of the main Theorem presented in Section 5. Section 6 considers additional permutations corresponding to non-simple Verma module multiplicities that are not obtained via the  construction described in Theorem~\ref{Main_Theorem}.

\section{Preliminaries}
\begin{small}
Throughout,  $\NN$ denotes the set of natural numbers, and $n\in \NN$.  For  $a,b\in \NN$ with $a\leq b$, we write $\left[a,b\right]:=\left\{k\in \NN \mid a\leq k\leq b\right\}$, 
 and for a finite set $S$ we let $|S|$ denote its cardinality. The symmetric group on $[1,n]$, denoted $S_n$, consists of all bijections on $[1,n]$ with composition $\cdot $ defined by $\left(\sigma\cdot  \tau\right)(k)=\sigma (\tau(k))$. The standard notation of a \textit{permutation} $\sigma\in S_n$  is $\sigma = \begin{scriptsize}\left(
\begin{array}{ccccc}
               1 &       2  & \cdots & n\\
        \sigma(1)&\sigma(2)& \cdots &\sigma(n)
\end{array}
\right)\end{scriptsize}$ for which we also use the
one-line notation $\sigma =\left(
        \sigma(1) \ \sigma(2) \ \cdots  \ \sigma(n)
\right)$. The \textit{inverse} of $\sigma$ is $\sigma^{-1}$. Throughout, $(i,j)$ denotes the \textit{transposition} interchanging $i$ and $j$, and we set $s_i := (i, i+1)$ for the simple transposition.
Every permutation admits an expression as a product of simple transpositions.  The \textit{length} of $\sigma \in S_{n}$, denoted by $l(\sigma)$, is the minimal integer $l$ such that 
$\sigma = s_{i_1}\cdot  s_{i_2} \cdots  s_{i_l}  $,  in which case  $ s_{i_1}\cdot  s_{i_2} \cdots  s_{i_l} $ is called a \emph{reduced decomposition} of 
$\sigma$. 
The relevant material can be found in \cite[Chapter 1]{Bjor}.  
  \end{small}

\begin{numeris}\normalfont{\textbf{Recall.}}
 We call  the pair $(i,j)$  an \textit{inversion} of $\sigma\in S_{n}$, if    $i<j$ and    $\sigma(i)>\sigma(j)$, or if    $i>j$ and    $\sigma(i)<\sigma(j)$. 
Let $\Inv(\sigma)$ be the set of inversions of $\sigma$, then 
\begin{center}
$l(\sigma) = \min\left\{l \mid  \sigma = s_{i_1}\cdot  s_{i_2} \cdots  s_{i_l}\right\} = \left|\Inv(\sigma)\right|$
\end{center}
 (see  \cite[Proposition 1.5.2]{Bjor}). For example, if  $\sigma =   \left(321456
\right)= s_1 \cdot  s_2\cdot  s_1\in S_6$, then  $\Inv(\sigma) =  \left\{(1,2),(1,3),(2,3)\right\}$, thus    $l(\sigma)=3$. 
The identity permutation $\textrm{id}_{n}:=(1\ 2 \ 3 \ \ldots \  n)$ has no inversions, thus $l(\textrm{id}_{n})=0$. For the  unique  longest element, $\omega_{n}:=\left(
\begin{array}{c}
        n\ n-1\  \cdots \ 2 \ 1
\end{array}
\right)$, in $S_{n}$ we have   $l(\omega_{n})=\left|\Inv(\omega_{n})\right|= \frac{n\cdot (n-1)}{2}$. 

The notation $(i,j)$ is used interchangeably to denote both, transpositions and inversions, depending on context. Since the distinction does not affect the structure of our arguments, this should not lead to confusion.
\label{inversions}
\end{numeris}

\begin{numeris}\normalfont{\textbf{Remark}.} 
The multiplication of $(i,j)$ with $\sigma$ from the left  interchanges  $i$ with $j$ (within the bottom line of $\sigma$) and from the right interchanges the positions of $\sigma(i)$ with $\sigma(j)$ and leaves everything else fixed.  
\begin{small}For  example $ (\textbf{2},\textbf{4})\cdot \left(
       \textbf{2}31 \textbf{4}
\right) = \left(
       \textbf{4}31 \textbf{2}
\right)$ and $ \left(
       2\textbf{3}1 \textbf{4}
\right)\cdot  (\textbf{\textit{2}},\textbf{\textit{4}}) =\left(
       2\textbf{4}1 \textbf{3}
\right)$. \end{small} Moreover 
\begin{center}
$\sigma \cdot  (i,j) = (\sigma(i), \sigma(j)) \cdot  \sigma$ \  \  \  as well as \  \  \  $(i,j) \cdot  \sigma  =  \sigma \cdot  (\sigma^{-1}(i), \sigma^{-1}(j)) $
\end{center}
\begin{small}For  example $\begin{scriptsize} (\textbf{2},\textbf{4})\cdot \left(
\begin{array}{ccccc}
       \textbf{2}&3&1& \textbf{4}
\end{array}
\right) = \left(
\begin{array}{ccccc}
       \textbf{2}&3&1& \textbf{4}
\end{array}
\right)\cdot  (\textbf{\textit{1}},\textbf{\textit{4}}) = \left(
\begin{array}{ccccc}
       4&3&1& 2
\end{array}
\right) \end{scriptsize}$ \end{small}.
\label{non_com} 
\end{numeris}

\begin{numeris}\normalfont{\textbf{Bruhat Order}.}
Let $(S,\leq)$ be a partially ordered set. We call $x\in S$ a \emph{smaller neighbor} of $y\in S$ if $x\neq y$ and, for every $z\in S$ with $x\leq z\leq y$, we have $x=z$ or $y=z$. (In the literature on posets, this relation is commonly referred to as a \emph{covering relation}.) The Bruhat order on $S_n$ will be defined in terms of such adjacent elements, for which we use the terminology \emph{smaller neighbor} (resp. \emph{larger neighbor}). For references see \cite[Chap. 2]{Bjor}, \cite[5.9-5.11]{Humph}, \cite[2.2]{Fed}.

\begin{subnum}\normalfont{\textbf{Definition}.} [\textbf{Bruhat order} on ($S_{n}, \ma$)]: Let $\sigma,\tau \in S_{n}$ with $\sigma \neq \tau$.
\begin{itemize}
	\item  We say that   $\sigma$ is a \textit{smaller neighbor}  of $\tau$, and  use the notation $\sigma \sma \tau $,  if  
	\begin{center}
	  $l(\tau) = l(\sigma) +1$   \  \    and  \  \  $\tau =  \sigma \cdot  (a,b)$   
	\end{center}
	 for some transposition  $(a,b)\in S_n$. Equivalently, $\tau$ \emph{covers} $\sigma$ in the Bruhat order.
	\item We define $\sigma < \tau$   if there exists a sequence $\rho_1, \ldots , \rho_t \in S_n $ \ such that  \ 
$\sigma = \rho_0\sma \cdots \sma \rho_t = \tau$
\end{itemize}
 \label{Bruhat}
\end{subnum}

\begin{small}For example $(13245)< (42315)$ in $S_5$, because $(13245)\sma \underbrace{(23415)}_{(13425)\cdot  (1,4)}\sma\underbrace{(32415)}_{(23415)\cdot  (1,2)}\sma \underbrace{(42315)}_{(3241)\cdot  (1,3)}$.  \end{small}

\begin{subnum}\normalfont{\textbf{Diagram.}}  Bruhat order can be visualized via diagram: the elements of $S_n$ are distributed according to their
 length (the length scale is on the right). Two elements connected by a line symbolize their neighborly status - the higher placed neighbor is always 
larger than the one below.  The following diagrams visually depict the Bruhat order in  $S_3$ and $S_4$.
  \label{Bruhat_Visualisation} 
\end{subnum}
\label{Bruhat_Order_all}
\end{numeris}
  
	\textrm{  }\\[-28mm]
  
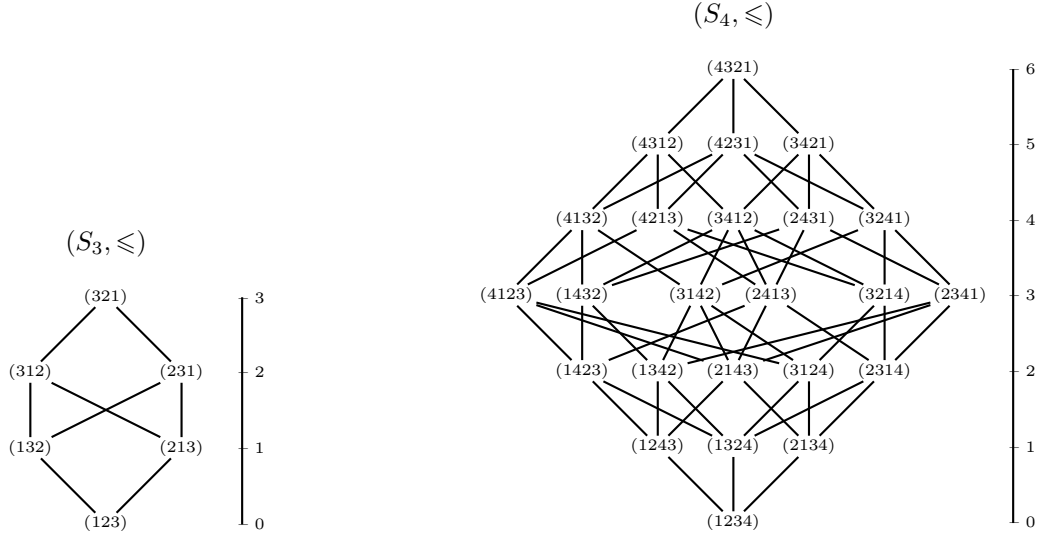
  \begin{tiny}
\psset{xunit=1mm,yunit=1mm,runit=1.5mm}
\begin{pspicture}(0,0)(0,0)
\rput(32,-50){
\rput(0,7){\begin{footnotesize}$\left(S_3,\ma \right)$\end{footnotesize}}
\rput(0,0){\rnode{123}{$(321)$}} 
\rput(-10,-10){\rnode{213}{$(312)$}} 
\rput(10,-10){\rnode{132}{$(231)$}}   
\rput(-10,-20){\rnode{231}{$(132)$}} 
\rput(10,-20){\rnode{312}{$(213)$}}
\rput(0,-30){\rnode{321}{$(123)$}}  
\rput(18,0.1){\rnode{oben3}{}}
\rput(18,-30.3){\rnode{unten3}{}}
\rput[l](17.5,-0.1){-- 3}
\rput[l](17.5,-10){-- 2}
\rput[l](17.5,-20){-- 1}
\rput[l](17.5,-30){-- 0}
}
             \psset{nodesep=0.5pt}
             \ncline{oben3}{unten3}
               \psset{nodesep=1pt,arrows=-}
 \psset{nrot=:U}              
\ncline{213}{123}  
\ncline{123}{132} 
\ncline{231}{213} 
\ncline{213}{312} 
\ncline{132}{231}  
\ncline{132}{312} 
\ncline{321}{231} 
\ncline{312}{321} 
\end{pspicture}
\\[17mm]  

\psset{xunit=1mm,yunit=1mm,runit=1.5mm}
\begin{pspicture}(0,0)(-115,-16)
\rput(0,7){\begin{footnotesize}$\left(S_4,\ma \right)$\end{footnotesize}}
              \rput(0,0){\rnode{1234}{$(4321)$}} 
              \rput(-10,-10){\rnode{2134}{$(4312)$}} 
              \rput(0,-10){\rnode{1324}{$(4231)$}}  
              \rput(10,-10){\rnode{1243}{$(3421)$}}  
               \rput(-20,-20){\rnode{2314}{$(4132)$}} 
               \rput(-10,-20){\rnode{3124}{$(4213)$}} 
               \rput(0,-20){\rnode{2143}{$(3412)$}}
               \rput(10,-20){\rnode{1342}{$(2431)$}}
               \rput(20,-20){\rnode{1423}{$(3241)$}}
               \rput(-30,-30){\rnode{3214}{$(4123)$}} 
               \rput(-20,-30){\rnode{2341}{$(1432)$}}
               \rput(-5,-30){\rnode{2413}{$(3142)$}}
               \rput(5,-30){\rnode{3142}{$(2413)$}}
               \rput(20,-30){\rnode{4123}{$(3214)$}}
               \rput(30,-30){\rnode{1432}{$(2341)$}}
               \rput(-20,-40){\rnode{3241}{$(1423)$}} 
               \rput(-10,-40){\rnode{2431}{$(1342)$}}
               \rput(0,-40){\rnode{3412}{$(2143)$}}
               \rput(10,-40){\rnode{4213}{$(3124)$}}
               \rput(20,-40){\rnode{4132}{$(2314)$}}
               \rput(-10,-50){\rnode{3421}{$(1243)$}} 
               \rput(0,-50){\rnode{4231}{$(1324)$}}
               \rput(10,-50){\rnode{4312}{$(2134)$}}
               \rput(0,-60){\rnode{4321}{$(1234)$}}  
               \rput(12,-31){\rnode{x}{}}
               \rput(-12,-31){\rnode{y}{}}
  \rput(37,0.1){\rnode{oben4}{}}
\rput(37,-60.3){\rnode{unten4}{}}
\rput[l](36.5,-0.1){-- 6}
\rput[l](36.5,-10){-- 5}
\rput[l](36.5,-20){-- 4}
\rput[l](36.5,-30){-- 3}
\rput[l](36.5,-40){-- 2}
\rput[l](36.5,-50){-- 1}
\rput[l](36.5,-60){-- 0}
\psset{nodesep=0.5pt}
\ncline{oben4}{unten4}
           \psset{nodesep=1pt,arrows=-}
\ncline{2134}{2314}
\ncline{1324}{3124}
\ncline{2143}{2341}
\ncline{1342}{3142}
\ncline{2413}{2431}
\ncline{1432}{3412}             
\ncline{4213}{4231}
\ncline{4132}{4312}
          \psset{nodesep=1pt,arrows=-}
\ncline{1243}{1423}
\ncline{1324}{1342}
\ncline{3124}{3142}
\ncline{2143}{4123}
\ncline{3214}{3412}
\ncline{2413}{4213}
\ncline{3241}{3421}
\ncline{2431}{4231}
    \psset{nodesep=1pt,arrows=-}
\ncline{2314}{2341} 
\ncline{2143}{2413} 
\ncline{1423}{4123} 
\ncline{3214}{3241} 
\ncline{3142}{3412} 
\ncline{1432}{4132} 

\psset{nodesep=1pt,arrows=-,linecolor=black}
\ncline{1234}{2134}
\ncline{1324}{2314}
\ncline{1243}{2143}
\ncline{3124}{3214}
\ncline{1342}{2341}
\ncline{1423}{2413}
\ncline{3142}{3241}
\ncline{4123}{4213}
\ncline{1432}{2431}
\ncline{3412}{3421}
\ncline{4132}{4231}
\ncline{4312}{4321}
               \psset{nodesep=1pt,arrows=-}
\ncline{1234}{1243}
\ncline{1324}{1423}
\ncline{2134}{2143}
\ncline{2314}{2413}
\ncline{3124}{4123}
\ncline{1342}{1432}
\ncline{3214}{4213}
\ncline{2341}{2431}
\ncline{3142}{4132}
\ncline{3241}{4231}
\ncline{3412}{4312}
\ncline{3421}{4321}
             \psset{nodesep=1pt,arrows=-}
\ncline{1234}{1324}
\ncline{2134}{3124}
\ncline{1243}{1342}
\ncline{2314}{3214}
\ncline{2143}{3142}
\ncline{1423}{1432}
\ncline{2341}{3241}
\ncline{2413}{3412}
\ncline{4123}{4132}
\ncline{2431}{3421}
\ncline{4213}{4312}
\ncline{4231}{4321}
\end{pspicture}
\end{tiny}
\\[4.5cm]

\section{\textbf{Critical Permutations}}
\begin{small}
In this section we introduce the class of \emph{critical} permutations, defined via the Bruhat order. We review known results for $S_n$ with $n \leq 5$ and present new computational data for $n=6,7$.  As noted in the introduction, these correspond to the simple modules in $\OO_{0}(\ger{sl}_{n}) $ where the principle Verma module has non-simple multiplicities.  
\end{small}

\begin{numeris}\normalfont{\textbf{Notation.}} Consider  $\tau \in S_{n}$ and the Bruhat order $\Bruh$ on $S_n$. We define
	\begin{center}
	$
		\begin{array}{ccl}
	\Lambda_{(\tau)} &:= & \left\{\sigma \in S_n  \mid  \sigma \Bruh \tau \right\}\\[2mm]
	\textbf{N}_{\textbf{s}}(\tau) &:=& \left\{\sigma \in S_{n} \mid \sigma \sma \tau \right\} 
	\\[3mm]
	\textbf{L}_{\textbf{t}}(\tau)&:=& \left\{\sigma \in S_{n} \mid \sigma \Bruh \tau \  \textrm{and} \ l(\sigma)=1 \right\} 
	\end{array}$
\end{center}
Thus, $\Lambda_{(\tau)}$ is the principal Bruhat ideal of $\tau$, $\mathbf{N}_{\mathbf{s}}(\tau)$ the set of its smaller neighbors, and $\mathbf{L}_{\mathbf{t}}(\tau)$ the set of simple transpositions comparable with $\tau$.
Since  $\textrm{id}_n\ma \sigma \ma \tau$ for each $\sigma \in \Lambda_{(\tau)}$ the permutation 
 $\tau$ is the unique maximal and $\textrm{id}_n$ the unique  minimal element in $\Lambda_{(\tau)}$ (visualized in the gray area in the following picture). 
\label{principal_set}
\end{numeris}

\begin{numeris}\normalfont{\textbf{Definition.}} We call  $\tau \in S_{n}$ a \textit{critical permutation of degree $\textbf{\textit{i}}$}, if

\begin{itemize}
	\item $\left|\left\{\eta \in \Lambda_{(\tau)}\mid l(\eta)=\textbf{\textit{i}} \right\}\right|\neq\left|\left\{\eta\in \Lambda_{(\tau)} \mid  l(\eta)=l(\tau)-\textbf{\textit{i}} \right\}\right|$ and 
	\item  $\left|\left\{\eta \in \Lambda_{(\tau)}\mid l(\eta)=\textbf{j} \right\}\right|=\left|\left\{\eta\in \Lambda_{(\tau)} \mid  l(\eta)=l(\tau)-\textbf{j} \right\}\right|$ for all $\textbf{j}\in \left[0,\textbf{\textit{i}}-1\right]$.
\end{itemize}

 We call $\tau \in S_{n}$ a \textit{critical permutation} if  $\tau $ is critical of degree $\textbf{\textit{i}}$ for some $\textbf{i}\in \left[0,  l(\tau)\right]$ 
 \label{Def_of_n_degree}
\end{numeris}

In other words: $\tau\in S_n $ is critical if for some  $0 \leq \textbf{i}\leq l(\tau)$  the number of  $\textbf{v}\in\Lambda_{(\tau)}$ with $l(\textbf{v})=\textbf{i}$ is \textit{not the same} as the number of  $\textbf{u}\in\Lambda_{(\tau)}$ with $l(\textbf{u})=l(\tau) - \textbf{i}$  (see the picture below).

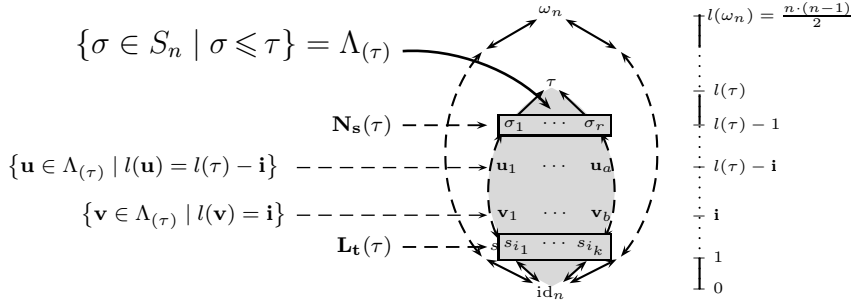
\textrm{  }\\[-1cm]

\hspace*{1.55cm}\psset{xunit=0.9cm,yunit=0.9cm,runit=1cm}
\begin{pspicture}(0,0)(0,0)
\rput(7.63,-2.89){
\psclip{
\pscircle[linestyle=none](-0.68,0){1.45cm}
\pscircle[linestyle=none](0.68,0){1.45cm}
        }
        \psframe*[linecolor=gray!30](-1,-1.47)(7,7)
\endpsclip
}
\end{pspicture}

\hspace*{1.55cm}\psset{xunit=0.7mm,yunit=0.7mm,runit=1mm}
\begin{pspicture}(0,30)(0,0)
\rput(100,17){
\begin{tiny}
\rput[l](28,9.6){\rnode{oben}{}}
\rput[l](28,2){\rnode{obenx1}{}}
\rput[l](28,-43.6){\rnode{unten}{}}
\rput[l](28,-36.4){\rnode{untenx1}{}}
\rput[l](26.5,9){\rnode{oben1}{$- l(\omega_{n})=\frac{n\cdot(n-1)}{2}$}}
\rput[l](26.5,-5.5){\rnode{oben3}{$- \ l(\tau)$}}
\rput[l](28,-5.5){\rnode{oben3x}{}}
\rput[l](26.5,-12){\rnode{oben4}{$- \ l(\tau)-1$}}
\rput[l](26.5,-20){\rnode{oben4i}{$- \ l(\tau)-\textbf{i}$}}
\rput[l](28,-12.5){\rnode{oben4x}{}}

\rput[l](26.5,-29){\rnode{oben4i}{$- \ \textbf{i}$}}
\rput[l](26.5,-37){\rnode{oben5}{$- \ 1 $}}
\rput[l](26.5,-43){\rnode{oben6}{$- \ 0 $}}
\rput(0,9.6){\rnode{A}{$\omega_n$}}
\rput(0,-43.6){\rnode{B}{$\textrm{id}_{n}$}}


\rput(20,-17){\rnode{aa}{$\textcolor{white}{a} $}}
\rput(-20,-17){\rnode{bb}{$\textcolor{white}{a} $}}
\rput(-13,2){\rnode{A1}{$$}}
\rput(-13,-37){\rnode{B1}{}}
\rput(13,2){\rnode{A2}{$$}}
\rput(13,-37){\rnode{B2}{}}

     
               \rput(0,-13.6){ \rput(0,9.6){\rnode{bA}{$\tau$}}
              \rput(0,-22){$\cdots$}
              \rput(0,1){$\cdots$}
              \rput(-9,-22){\rnode{bxx}{$s_{i_1}$}}
              \rput(9,-22){\rnode{bxy}{$s_{i_k}$}}
              \rput(-9,1){\rnode{b11}{}}
\rput(9,1){\rnode{b21}{}}   
    }
              \end{tiny}
                                    }
                      \psset{nodesep=1.3pt}
                 \ncline[linestyle=dotted]{-}{oben}{unten}   
                 \ncline{-}{oben}{obenx1}   
                 \ncline{-}{unten}{untenx1}  
                 \ncline{-}{oben3x}{oben4x}  
                     
              \ncline{<->}{3}{20}
             \ncline{<->}{3}{21}
             \ncline{<->}{3}{10}
             \ncline{<->}{3}{11}
              \ncline{<->}{A}{A1}
               \ncline{<->}{A}{A11}
                \ncline{<->}{A}{A21}
              \ncline{<->}{B1}{B}
              \ncline{<->}{A}{A2}
              \ncline{<->}{B2}{B}
              
               \ncline{<->}{bxx}{B}
                \ncline{<->}{bxy}{B}
               \ncline{<->}{bA}{b11}
                \ncline{<->}{bA}{b21}
              
              \ncarc[linewidth=0.5pt,arcangle=-15,linestyle=dashed]{<->}{sn}{20}
              \ncarc[linewidth=0.5pt,arcangle=15,linestyle=dashed]{<->}{ln}{10}
              \ncarc[arcangle=25,linestyle=dashed]{<->}{A21}{10}
              \ncarc[arcangle=-25,linestyle=dashed]{<->}{A11}{11}
              
              \ncarc[arcangle=25,linestyle=dashed]{<->}{b21}{bxy}
              \ncarc[arcangle=-25,linestyle=dashed]{<->}{b11}{bxx}
              
              \ncarc[arcangle=-25,linestyle=dashed]{<->}{B}{20}
              \ncarc[arcangle=25,linestyle=dashed]{<->}{B}{21}
               \psset{nodesep=0.3pt}
              \ncarc[arcangle=20,linestyle=dashed]{<-}{A2}{aa}
              \ncarc[arcangle=20,linestyle=dashed]{->}{aa}{B2}
              \ncarc[arcangle=-20,linestyle=dashed]{<-}{A1}{bb}
              \ncarc[arcangle=20,linestyle=dashed]{<-}{B1}{bb}
\end{pspicture} 
\psset{xunit=0.7mm,yunit=0.7mm,runit=1mm}
\begin{pspicture}(0,0)(0,0)
\rput(97,11){
              \rput(-28,8){\rnode{v2}{}}
              \rput(0,-3){\rnode{v3}{}}
              \rput[r](-30,9){\rnode{v4}{$\left\{\sigma \in S_{n}\mid \sigma \Bruh \tau\right\}=\Lambda_{(\tau)}$}}
 \ncarc[linewidth=1pt,arcangle=25,linecolor=black]{->}{v2}{v3}
\rput[r](-30,-6){\begin{scriptsize}$\textbf{N}_{\textbf{s}}(\tau)$\end{scriptsize}}
							\psline[linewidth=0.7pt,linestyle=dashed]{->}(-28,-6)(-12,-6)
\rput[r](-30,-29){\begin{scriptsize}$\textbf{L}_{\textbf{t}}(\tau)$\end{scriptsize}}
							\psline[linewidth=0.7pt,linestyle=dashed]{->}(-28,-29)(-12,-29)

 \begin{tiny} \rput(0.7,-23){\rnode{dez-i}{$\textbf{v}_1 \  \  \       \cdots \    \  \     \textbf{v}_b  $}}
\rput(0.7,-6){\rnode{dez-ill}{\setlength{\fboxsep}{1.7pt}{\fcolorbox{black}{gray!30}{$\sigma_1 \   \      \cdots \    \    \sigma_r  $}}}}
             \rput(0.7,-14){\rnode{dez-m-i}{$\textbf{u}_1 \  \  \       \cdots \    \  \     \textbf{u}_a  $}} 
						\rput(0.7,-29){\rnode{dez-m-i}{\setlength{\fboxsep}{1.7pt}{\fcolorbox{black}{gray!30}{$s_{i_1} \         \cdots \       s_{i_k} $}}}}
						\end{tiny}
	\psline[linewidth=0.5pt,linestyle=dashed]{->}(-48,-23)(-12,-23)
	\psline[linewidth=0.5pt,linestyle=dashed]{->}(-48,-14)(-12,-14)
	\rput[r](-50,-14){\begin{scriptsize}$\left\{\textbf{u} \in \Lambda_{(\tau)} \mid l(\textbf{u}) =l( \tau)-\textbf{i}\right\}$
	\end{scriptsize}}
	\rput[r](-50,-23){\begin{scriptsize}$\left\{\textbf{v} \in \Lambda_{(\tau)} \mid l(\textbf{v}) =\textbf{i}\right\}$\end{scriptsize}}
}
\end{pspicture} 

\textrm{  }\\[1.9cm]	
It should be noted that no link is currently known between critical permutations of degree $i$ and a more precise understanding of Verma multiplicities (see ~\ref{what_we_absurved_until_now} for further details).

The primary focus of this paper is on permutations of degree $\textbf{\textit{1}}$. We shall therefore refer to them as \textit{first-degree critical permutations} or simply \textit{critical permutations in the first-degree}, unless stated otherwise. 
By definition, $\tau\in S_n$ is first-degree critical permutation  if and only if

\begin{center}
$\left|\textbf{N}_{\textbf{s}}(\tau)\right| \neq \left|\textbf{L}_{\textbf{t}}(\tau) \right| $
\end{center}
Since $S_n$ contains exactly $n-1$ simple transpositions, $s_1, s_2, \ldots, s_{n-1}$, and $\textbf{L}_{\textbf{t}}(\tau) \subseteq S_n$, we always have $\left|\textbf{L}_{\textbf{t}}(\tau) \right| \leq n-1$.
This observation immediately yields the following result.

\begin{numeris}\normalfont{\textbf{Lemma.}} \textit{Let $\tau\in S_n$ and $\left|\textbf{N}_{\textbf{s}}(\tau)\right|\geq n$, then $\tau$ is a first-degree critical permutation.} \\

\textit{Proof.}  Since $\left|\textbf{N}_{\textbf{s}}(\tau)\right|\geq n$ and  $\left|\textbf{L}_{\textbf{t}}(\tau) \right|\leq n-1$, we have $\left|\textbf{N}_{\textbf{s}(\tau)}\right|\neq \left|\textbf{L}_{\textbf{t}}(\tau) \right| $. \hfill $\Box$ 
\label{firs_degree_for_sure}
\end{numeris}

\begin{numeris}\normalfont{\textbf{Examples of critical permutations.}} 
It is evident that  for   $\textrm{id}_n$, all simple transpositions, $s_i$,  and all permutation of length two, $s_{i_1}\cdot  s_{i_2}$ with $i_1\neq i_2$,
we have 
$\Lambda_{(\textrm{id}_n)} = \left\{\textrm{id}_n\right\}$,  $\Lambda_{(s_{i})} = \left\{  \textrm{id}_n, s_{i}\right\}$ and $\Lambda_{(s_{i_1}\cdot  s_{i_2})} = \left\{\textrm{id}_n, s_{i_1}, s_{i_2}, s_{i_1}\cdot  s_{i_2}\right\}$. For $\tau\in S_n$ with $l(\tau)=3$ we have  either  $\tau = s_i\cdot  s_{i\pm 1}\cdot  s_i$ or  $\tau = s_{i_1}\cdot  s_{i_2}\cdot  s_{i_3}$ where $i_1, i_2, i_3$ are pairwise different. Thus  $\Lambda_{(s_i\cdot  s_{i\pm 1}\cdot  s_i)}$ has 6 elements and $\Lambda_{(s_{i_1}\cdot  s_{i_2}\cdot  s_{i_3})}$ has 8 elements represented in the Bruhat-diagrams below. The "up-down"-symmetry of the diagrams reveal that there  are no  critical permutations with  $ l(\tau)\leq 3$. \\[2.5mm]

\begin{tiny}
\psset{xunit=1mm,yunit=1mm,runit=1.5mm}
\begin{pspicture}(0,0)(0,0)
\rput(2,-5){
\rput(0,0){
 \rput(0,-10){\rnode{ssi}{$s_{i}$}}
\rput(0,-20){\rnode{idi}{$\textrm{id}_n$}}
}
\rput(30,0){
\rput(0,0){\rnode{s1s2}{$s_{i_1}\cdot  s_{i_2}$}}
\rput(-10,-10){\rnode{s1}{$s_{i_1}$}}  \rput(10,-10){\rnode{s2}{$s_{i_2}$}}
\rput(0,-20){\rnode{id}{$\textrm{id}_n$}}
}

\rput(70,0){
\rput(0,10){\rnode{s_121}{$s_{i}\cdot  s_{i\pm 1}\cdot  s_{i}$}}
\rput(-10,0){\rnode{s_12}{$s_{i}\cdot  s_{i\pm 1}$}} \rput(10,0){\rnode{s_21}{$s_{i\pm 1}\cdot  s_{i}$}}
\rput(-10,-10){\rnode{s_1}{$s_{i_1}$}}  \rput(10,-10){\rnode{s_2}{$s_{i_2}$}}
\rput(0,-20){\rnode{id_12}{$\textrm{id}_n$}}
}

\rput(115,0){
\rput(0,10){\rnode{c_123}{$s_{i_1}\cdot  s_{i_2}\cdot  s_{i_3}$}}
\rput(-15,0){\rnode{c_12}{$s_{i_1}\cdot  s_{i_2}$}} \rput(0,0){\rnode{c_13}{$s_{i_1}\cdot  s_{i_3}$}} \rput(15,0){\rnode{c_23}{$s_{i_2}\cdot  s_{i_3}$}}
\rput(-15,-10){\rnode{c_1}{$s_{i_1}$}}  \rput(0,-10){\rnode{c_2}{$s_{i_2}$}} \rput(15,-10){\rnode{c_3}{$s_{i_3}$}}
\rput(0,-20){\rnode{id_123}{$\textrm{id}_n$}}
}

\rput(150,0){ 
\psline[linewidth=0.7pt]{-}(0,11)(0,-21)
              \rput[l](-1.2,10){\begin{scriptsize}$- \ \ 3 $\end{scriptsize}}
							\rput[l](-1.2,0){\begin{scriptsize}$- \ \ 2 $\end{scriptsize}}
							\rput[l](-1.2,-10){\begin{scriptsize}$- \ \ 1 $\end{scriptsize}}
							\rput[l](-1.2,-20){\begin{scriptsize}$- \ \ 0 $\end{scriptsize}}}
}
\psset{nodesep=1pt,arrows=-}
\ncline{ssi}{idi}
\ncline{s1s2}{s1}\ncline{s1s2}{s2} \ncline{id}{s1} \ncline{id}{s2}
\ncline{id_12}{s_2} \ncline{id_12}{s_1} \ncline{s_12}{s_2} \ncline{s_12}{s_1} \ncline{s_21}{s_2} \ncline{s_21}{s_1} \ncline{s_12}{s_121} \ncline{s_21}{s_121}
\ncline{id_123}{c_1} \ncline{id_123}{c_2} \ncline{id_123}{c_3} \ncline{c_1}{c_12} \ncline{c_1}{c_13} \ncline{c_2}{c_12} \ncline{c_2}{c_23} \ncline{c_3}{c_13} \ncline{c_3}{c_23} \ncline{c_12}{c_123} \ncline{c_13}{c_123}  \ncline{c_23}{c_123} 

\end{pspicture}\end{tiny}\\[2.1cm] 

\begin{subnum}\normalfont{\textbf{There are no critical permutations in $S_n$ for $n=1,2,3$.}} Since $ l(\tau)\leq 3$ for all $\tau \in S_n$ with $1\leq n \leq 3$, there are no critical permutations in $S_1$, $S_2$, and $S_3$.
 \label{no_critical_in_S_1_3}
\end{subnum}

\begin{subnum}\normalfont{\textbf{There are \textit{two} critical permutations in $S_4$.}} 
  The first time  critical permutations appear in $S_4$:
There are exactly two permutation, $\textbf{x}_1 = (3412)$ and $\textbf{x}_2=(4231)$, in $S_{4}$ that are critical. Since
 $\underbrace{\left|\left\{\sigma\in \Lambda_{(\textbf{x}_j)} \mid  l(\sigma)=l(\textbf{x}_j)-\textbf{1} \right\}\right|}_{4} \neq \underbrace{\left|\left\{\sigma\in \Lambda_{(\textbf{x}_j)} \mid  l(\sigma)=\textbf{1} \right\}\right|}_{3} $,  both permutations are critical in first-degree. The diagrams of $(\Lambda_{(\textbf{x}_j)}, \Bruh)$ are represented   below; here $j=1,2$.
 \label{crit_in_S_4}
\end{subnum}

\textrm{  }
\\[3mm]

\begin{tiny}
\psset{xunit=1mm,yunit=1mm,runit=1.5mm}
\begin{pspicture}(0,0)(0,0)

\rput(20,20){ \rput(-25,-32.2){\psframe*[linecolor=gray!30](0,0)(50,4.4)}   \rput(-25,-32.2){\psframe[linecolor=black!50](0,0)(50,4.4)} 
             \rput(-15,-52.2){\psframe*[linecolor=gray!30](0,0)(30,4.4)}    \rput(-15,-52.2){\psframe[linecolor=black!50](0,0)(30,4.4)}     
               \rput(0,-20){\rnode{2143}{$(3412)$}}
               \rput(-20,-30){\rnode{2341}{$(1432)$}}
               \rput(-5,-30){\rnode{2413}{$(3142)$}}
               \rput(5,-30){\rnode{3142}{$(2413)$}}
               \rput(20,-30){\rnode{4123}{$(3214)$}}   
               \rput(-20,-40){\rnode{3241}{$(1423)$}}
               \rput(-10,-40){\rnode{2431}{$(1342)$}}
               \rput(0,-40){\rnode{3412}{$(2143)$}}
               \rput(10,-40){\rnode{4213}{$(3124)$}}
               \rput(20,-40){\rnode{4132}{$(2314)$}}
               \rput(-10,-50){\rnode{3421}{$(1243)$}}
               \rput(0,-50){\rnode{4231}{$(1324)$}}
               \rput(10,-50){\rnode{4312}{$(2134)$}}
               \rput(0,-60){\rnode{4321}{$(1234)$}}
               \rput(12,-31){\rnode{x}{}}
               \rput(-12,-31){\rnode{y}{}}
							\rput(-4,0){\psline[linewidth=0.7pt]{-}(35,-19)(35,-61)
							\rput[l](33.8,-20){\begin{scriptsize}$- \ \ 4= l(\textbf{x}_{1}) $\end{scriptsize}}
							\rput[l](33.8,-30){\begin{scriptsize}$- \ \ 3 $\end{scriptsize}}
							\rput[l](33.8,-40){\begin{scriptsize}$- \ \ 2 $\end{scriptsize}}
							\rput[l](33.8,-50){\begin{scriptsize}$- \ \ 1 $\end{scriptsize}}
							\rput[l](33.8,-60){\begin{scriptsize}$- \ \ 0 $\end{scriptsize}}}
						}
\psset{nodesep=0.5pt}
\ncline{oben4}{unten4}
           \psset{nodesep=2pt,arrows=-}
\ncline{2143}{2341}
\ncline{2413}{2431}
\ncline{4213}{4231}
\ncline{4132}{4312}
          \psset{nodesep=2pt,arrows=-}
\ncline{2143}{4123}
\ncline{2413}{4213}
\ncline{3241}{3421}
\ncline{2431}{4231}
    \psset{nodesep=2pt,arrows=-}
\ncline{2143}{2413} 
\ncline{3142}{3412}  

\psset{nodesep=2pt,arrows=-,linecolor=black}
\ncline{3142}{3241}
\ncline{4123}{4213}
\ncline{3412}{3421}
\ncline{4132}{4231}
\ncline{4312}{4321}
               \psset{nodesep=2pt,arrows=-}
\ncline{2341}{2431}
\ncline{3142}{4132}
\ncline{3241}{4231}
\ncline{3412}{4312}
\ncline{3421}{4321}
             \psset{nodesep=2pt,arrows=-}
\ncline{2143}{3142}
\ncline{2341}{3241}
\ncline{2413}{3412}
\ncline{4123}{4132}
\ncline{2431}{3421}
\ncline{4213}{4312}
\ncline{4231}{4321}


\rput(110,20){
\rput(-25,-22.2){\psframe*[linecolor=gray!30](0,0)(50,4.4)}   \rput(-25,-22.2){\psframe[linecolor=black!50](0,0)(50,4.4)} 
\rput(-15,-52.2){\psframe*[linecolor=gray!30](0,0)(30,4.4)}    \rput(-15,-52.2){\psframe[linecolor=black!50](0,0)(30,4.4)}   
\rput(0,-10){\rnode{4231a}{$(4231)$}}
               \rput(-20,-20){\rnode{4132a}{$(4132)$}} \rput(-5,-20){\rnode{4213a}{$(4213)$}} \rput(5,-20){\rnode{2431a}{$(2431)$}} \rput(20,-20){\rnode{3241a}{$(3241)$}}
							
							\rput(-30,-30){\rnode{4123a}{$(4123)$}} \rput(-15,-30){\rnode{1432a}{$(1432)$}} \rput(-5,-30){\rnode{3142a}{$(3142)$}} \rput(5,-30){\rnode{2413a}{$(2413)$}} \rput(15,-30){\rnode{3214a}{$(3214)$}} \rput(30,-30){\rnode{2341a}{$(2341)$}}

              \rput(-20,-40){\rnode{1423a}{$(1423)$}}
               \rput(-10,-40){\rnode{1342a}{$(1342)$}}
               \rput(0,-40){\rnode{2143a}{$(2143)$}}
               \rput(10,-40){\rnode{3124a}{$(3124)$}}
               \rput(20,-40){\rnode{2314a}{$(2314)$}}

               \rput(-10,-50){\rnode{1243a}{$(1243)$}}
               \rput(0,-50){\rnode{1324a}{$(1324)$}}
               \rput(10,-50){\rnode{2134a}{$(2134)$}}
               \rput(0,-60){\rnode{1234a}{$(1234)$}}
               \rput(12,-31){\rnode{x}{}}
               \rput(-12,-31){\rnode{y}{}}
							\rput(-4,0){\psline[linewidth=0.7pt]{-}(40,-9)(40,-61)
							\rput[l](38.8,-10){\begin{scriptsize}$- \ \ 5= l(\textbf{x}_{2}) $\end{scriptsize}}
							\rput[l](38.8,-20){\begin{scriptsize}$- \ \ 4$\end{scriptsize}}
							\rput[l](38.8,-30){\begin{scriptsize}$- \ \ 3 $\end{scriptsize}}
							\rput[l](38.8,-40){\begin{scriptsize}$- \ \ 2 $\end{scriptsize}}
							\rput[l](38.8,-50){\begin{scriptsize}$- \ \ 1 $\end{scriptsize}}
							\rput[l](38.8,-60){\begin{scriptsize}$- \ \ 0 $\end{scriptsize}}}
						}

\psset{nodesep=0.5pt}
\ncline{oben4}{unten4}
           \psset{nodesep=2pt,arrows=-}
\ncline{1234a}{2134a} \ncline{1234a}{1324a} \ncline{1234a}{1243a}
\ncline{2134a}{2143a} \ncline{2134a}{3124a} \ncline{2134a}{2314a}
\ncline{1324a}{1423a} \ncline{1324a}{1342a} \ncline{1324a}{3124a} \ncline{1324a}{2314a}
\ncline{1243a}{1423a} \ncline{1243a}{1342a} \ncline{1243a}{2143a}

\ncline{1423a}{4123a} \ncline{1423a}{1432a} \ncline{1423a}{2413a}
\ncline{1342a}{1432a} \ncline{1342a}{3142a} \ncline{1342a}{2341a}
\ncline{2143a}{4123a} \ncline{2143a}{3142a} \ncline{2143a}{2413a} \ncline{2143a}{2341a} 
\ncline{3124a}{4123a} \ncline{3124a}{3142a} \ncline{3124a}{3214a} 
\ncline{2314a}{2413a} \ncline{2314a}{3214a} \ncline{2314a}{2341a}
\ncline{4231a}{4132a} \ncline{4231a}{4213a} \ncline{4231a}{2431a} \ncline{4231a}{3241a}

\ncline{4132a}{4123a}  \ncline{4132a}{1432a}  \ncline{4132a}{3142a} 
\ncline{4213a}{4123a}  \ncline{4213a}{2413a}  \ncline{4213a}{3214a} 
\ncline{2431a}{1432a}  \ncline{2431a}{2413a}  \ncline{2431a}{2341a}  
\ncline{3241a}{3142a} \ncline{3241a}{2314a} \ncline{3241a}{2341a}

\end{pspicture}

\end{tiny}
\textrm{  }\\[3.7cm]

\begin{subnum}\normalfont{\textbf{Critical permutations in $S_5$.}} In \cite[5.24]{Jantzen} Jantzen  identifies $32$ elements $q \in S_5$ with multiplicities $[M(\textrm{id}_5),L(q)]\neq 1$. These correspond precisely to the critical permutations of $S_5$, which are displayed in the symbolic Bruhat diagram in Section~\ref{critical_26}.
 With the exception of $(45312)$, that is critical in \textit{second} degree,  all remaining 31 elements are  critical in first degree. 
 \label{32_in_S_5}
\end{subnum}

\begin{subnum}\normalfont{\textbf{Critical permutations in $S_6$ and in $S_7$.}} 
Based on the computational program developed by Elbosari (see \cite{Elbos}), there are $354$ critical permutations in $S_6$ and $3,488$ in $S_7$. The relationship between these permutations, together with further details, is explored in ~\ref{Critical_in_S_7}. 
 \label{critical_in_S_6}
\end{subnum}

It is remarkable  that for $n=4,5,6,7$ the permutation $q_{n-3}\in S_{n}$ of the form

\begin{center}
$q_{n-3}:=\begin{scriptsize} \left(\begin{array}{cc}
	1 &2 \\
	n-1 & n 
\end{array}
\framebox{$
\begin{array}{ccccc}
	 3 & 4 & \cdots & n-3 & n-2 \\
	 n-2 & n-3& \cdots &  4 &  3  
\end{array}$}\begin{array}{cc}
	n-1 & n \\
	 1 & 2 
\end{array}
\right)\end{scriptsize}$ 
\end{center}
 is   critical  of degree $n-3$, and for $n=5,6,7$  the only critical permutation in this degree. This observation suggests the possibility that a critical permutation $q_{\textit{\textbf{i}}}$ is  of degree \textit{\textbf{i}}  for all $\textit{\textbf{i}}\geq 2$.
\label{critical_example}
\end{numeris}

\section{Main Theorem}
\begin{small}
In this section we present the relationship between the Bruhat order of $S_n$ and $S_{n+1}$, which establishes a link between the   \textit{critical} permutations of $S_{n}$ and  $S_{n+1}$. Using this relationship, we observe that $26$ of the $32$ critical permutations in $S_5$ are derived from  $(3412)$ and $(4231)$ in $S_4$ (see ~\ref{critical_example}). Continuing this procedure, the critical permutations in $S_5$ generate $343$ of the $345$  permutations in $S_6$, while only $20$ of the $3,488$ critical permutations in $S_7$ cannot be obtained from those in  $S_6$.
\end{small}

\begin{numeris}\normalfont{\textbf{Cosets} $\textbf{I}_{(b)}^{(a)}$ \textbf{of} $S_{n+1}$}.  It is well known that   $(S_{n}, \cdot )$ is a subgroup of $(S_{n+1}, \cdot )$, so  $\sigma\in S_{n}$ can be considered  as a permutation in  $S_{n+1}$ with $\sigma(n+1) = n+1$. Moreover, $S_{n+1}$ can be expressed as a union of right cosets with respect to $S_n$ as well as left cosets  (see \cite[Chapter 4,5]{Humphreys}). For   $a,b\in \left[1,n+1\right]$  denote by $\textbf{I}_{(b)}^{(a)}$ the set of permutations $\sigma\in S_{n+1}$ with $\sigma(a)=b$.  Obviously  $\textbf{I}_{(b)}^{(a)}\cap \textbf{I}_{(b')}^{(a)} = \emptyset$  for $b\neq b'$,  and  $\textbf{I}_{(b)}^{(a)}\cap \textbf{I}_{(b)}^{(a')} = \emptyset$ for $a\neq a'$.  
For  $i\in \left[1,n\right]$ denote

\begin{center}
$c_{i}:= s_{i}\cdot  s_{i+1}\cdots  s_{n-1} \cdot  s_{n} $  
 \  \  \  \  \  and  \  \  \   \  \   $d_{i}:=s_{n}\cdot  s_{n-1}\cdots  s_{i+1}\cdot  s_{i} $ 
\end{center}
where $s_{j}=(j,j+1)$ denotes the adjacent transposition. Additionally, set $c_{n+1}=d_{n+1}=\textrm{id}_{n+1}$.
 Note that  $c_{i}= s_{i}\cdot  c_{i+1}$, and $(c_{i})^{-1} = d_i$, thus $d_{i}= d_{i+1}\cdot  s_{i}$.  Moreover,   $(c_{i}\cdot  \sigma)(n+1)=i $   and     $(\sigma \cdot  d_{i})(i)=n+1 $ for each  $\sigma\in S_n=\textbf{I}_{(n+1)}^{(n+1)}$. The bijective functions $S_n \xrightarrow{c_{i}\cdot  -} \textbf{I}_{(i)}^{(n+1)}$,  $S_n \xrightarrow{ -\cdot  d_{i}} \textbf{I}^{(i)}_{(n+1)}$,  
 $\textbf{I}_{(i+1)}^{(n+1)}\xrightarrow{s_{i}\cdot  -} \textbf{I}_{(i)}^{(n+1)}$,  $\textbf{I}^{(i+1)}_{(n+1)}\xrightarrow{ -\cdot  s_{i}} \textbf{I}^{(i)}_{(n+1)}$ and   $\textbf{I}_{(i)}^{(n+1)} \xrightarrow{ (-)^{-1}}  \textbf{I}^{(i)}_{(n+1)}$ imply

\begin{center}
$
\begin{array}{cclclcl}
\textbf{I}_{(i)}^{(n+1)} &=&c_i\cdot  S_n = \left\{ c_{i}\cdot  \sigma  \mid \sigma \in S_{n} \right\} &=&  \left\{\eta\in S_{n+1} \mid \eta(n+1)=i\right\}  
&=& \left\{s_{i}\cdot  \eta  \mid \eta \in \textbf{I}_{(i+1)}^{(n+1)}\right\} \\[3mm]
\textbf{I}^{(i)}_{(n+1)} & = &  S_n  \cdot  d_i = \left\{\eta^{-1} \mid \eta \in \textbf{I}_{(i)}^{(n+1)}\right\}
&=& \left\{\eta\in S_{n+1} \mid \eta(i)=n+1\right\} 
&=& \left\{\eta \cdot  s_{i} \mid \eta \in \textbf{I}^{(i+1)}_{(n+1)}\right\}.
\end{array}$   
\end{center}
Thus 
	$S_{n+1} = \displaystyle \bigcup_{i=1}^{n+1}  \ \textbf{I}^{(n+1)}_{(i)} =  \displaystyle \bigcup_{i=1}^{n+1} \  \textbf{I}^{(i)}_{(n+1)}$
	allows us to view  $S_{n+1}$  as $n+1$ "\textit{copies}" of $S_n$. In  ~\ref{Bruhat_order_on_Blocks}  we establish  that  for each $\textbf{x}\in \textbf{I}_{(i)}^{(n+1)} $ we have $\textbf{x}\sma \left(s_{i-1}\cdot \textbf{x}\right)$, which in turn implies $l(s_{i-1}\cdot \textbf{x})= l(\textbf{x})+1$. Therefore the permutations in \textit{block} $\textbf{I}_{(i-1)}^{(n+1)}$ in the Bruhat-order diagram of $S_{n+1}$ are obtained from those in $\textbf{I}_{(i)}^{(n+1)}$ by a one-step shift on the length scale, as represented in the diagram below. Moreover, if $\sigma\sma \tau $  for $\sigma, \tau\in S_n$, then for the corresponding permutations in $\textbf{I}_{(i)}^{(n+1)} $ we have  $\left(s_{i}\cdot \sigma\right)\sma \left(s_{i}\cdot \tau\right) $ as well, which implies that the Bruhat-order "neighborhood relationships" within the block $\textbf{I}_{(i)}^{(n+1)}$ follow the same pattern as in $S_n$ (see  Lemma~\ref{Bruhat_order_on_Blocks}(2)).

  Notice that the  diagram below does not accurately represent the Bruhat-order of $S_{n+1}$, as it omits the \textit{neighboring-lines} between the elements in  $\textbf{I}_{(i)}^{(n+1)}$ and $\textbf{I}_{(j)}^{(n+1)}$ where $\left|i-j\right|\geq 2$.

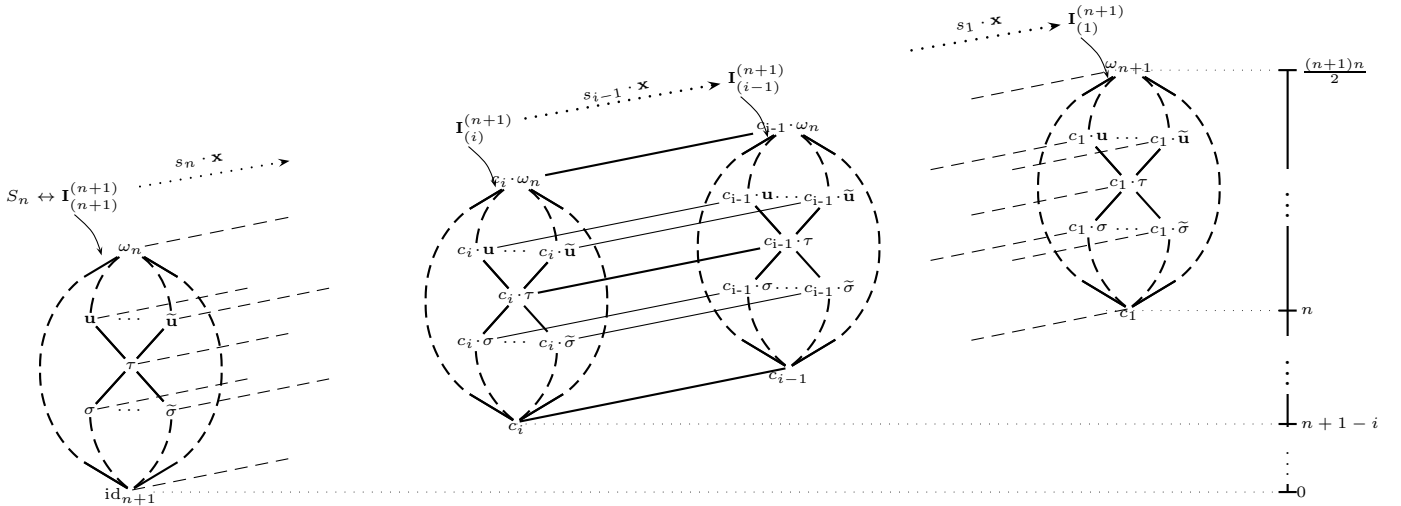
\textrm{  }\\

\psset{xunit=0.6mm,yunit=0.6mm,runit=1mm}
\begin{pspicture}(0,5)(0,0)
\rput(55,-48){
\rput(220,-49){
\psline[linestyle=solid,linecolor=black]{-}(0,0)(0,93) \psline[linestyle=solid,linecolor=black]{-}(-2,93)(2,93) \rput[l](3,93){\begin{tiny}$\frac{(n+1)n}{2}$\end{tiny}} \psline[linestyle=dotted,linewidth=0.4pt](-40,93)(0,93)
\psline[linestyle=solid,linecolor=black]{-}(-2,0)(2,0) \rput(3,0){\begin{tiny}0\end{tiny}}\psline[linestyle=dotted,linewidth=0.4pt](-170,15)(0,15)\psline[linestyle=solid,linecolor=black]{-}(-2,15)(2,15) \rput[l](3,15){\begin{tiny}$n+1-i$\end{tiny}}  \psline[linestyle=solid,linecolor=black]{-}(-2,40)(2,40) \psline[linestyle=dotted,linewidth=0.4pt](-40,40)(0,40)    \rput[l](3,40){\begin{tiny}$n$\end{tiny}}
\psline[linestyle=dotted,linewidth=0.4pt]{-}(-250,0)(0,0)
\rput(0,65){\colorbox{white}{$\vdots$}}  \rput(0,28){\colorbox{white}{$\vdots$}} \rput(0,8){\begin{tiny}\colorbox{white}{$\vdots$}\end{tiny}}
}
\rput(-35,-5){\rput(0,23){\psline[linestyle=dotted]{->}(0,0)(35,6)\rput{8}(15,6){\begin{tiny}$s_n\cdot  \textbf{x}$\end{tiny}}
}
\begin{tiny}
\rput(-15,21){$S_n \leftrightarrow \textbf{I}_{(n+1)}^{(n+1)}$\pscurve[linewidth=0.4pt,linestyle=solid]{->}(-10,-3)(-6,-7)(-4,-12)}
\rput(0,9.6){\rnode{A}{$\omega_{n}$}}
\rput(0,-43.6){\rnode{B}{}}
\rput(0,-45){$\textrm{id}_{n+1}$}
\rput(20,-17){\rnode{aa}{}}
\rput(-20,-17){\rnode{bb}{}}
\rput(-11,3){\rnode{A1}{}}
\rput(-11,-37){\rnode{B1}{}}
\rput(11,3){\rnode{Ax}{}}
\rput(11,-37){\rnode{Bx}{}} 
\rput(35,7){\rput(0,9.6){\rnode{An}{}}
\rput(0,-43.6){\rnode{Bn}{}}
\rput(9,-26){\rnode{20n}{}}
              \rput(-9,-26){\rnode{21n}{}}
              \rput(0,-16){\rnode{3n}{}}      
              \rput(9,-6){\rnode{10n}{}}
              \rput(-9,-6){\rnode{11n}{}}
}
              
              \rput(9,-26){\rnode{20}{$\widetilde{\sigma}$}}
              \rput(0,-26){$\cdots$}
              \rput(-9,-26){\rnode{21}{$\sigma$}}
              \rput(0,-16){\rnode{3}{$\tau$}}      
              \rput(9,-6){\rnode{10}{$\widetilde{\textbf{u}}$}}
              \rput(0,-6){$\cdots$}
              \rput(-9,-6){\rnode{11}{$\textbf{u}$}}
              \end{tiny}}
							
							\psset{nodesep=0.5pt,linestyle=dashed,linewidth=0.4pt}  
								\ncline{-}{A}{An}
								\ncline{-}{B}{Bn}
								\ncline{-}{20n}{20}
								\ncline{-}{21n}{21}
								\ncline{-}{3n}{3}
								\ncline{-}{10n}{10}
								\ncline{-}{11n}{11}
\rput(50,10){
\rput(0,21){\psline[linestyle=dotted,linewidth=1pt]{->}(0,2)(45,10)\rput{11}(22,8){\begin{tiny}$s_{i-1}\cdot  \textbf{x}$\end{tiny}}
}
\begin{tiny}
\rput(-7,21){$ \displaystyle \textbf{I}_{(i)}^{(n+1)}$\pscurve[linewidth=0.4pt,linestyle=solid]{->}(-10,-3)(-6,-7)(-4,-12)}  
\rput(0,9.6){\rnode{A2}{$c_{i}\!\cdot\!\omega_{n}$}}

\rput(0,-43.6){\rnode{B2}{}}
\rput(0,-45){$c_i$}
\rput(20,-17){\rnode{aa2}{}}
\rput(-20,-17){\rnode{bb2}{}}
\rput(-11,3){\rnode{A12}{}}
\rput(-11,-37){\rnode{B12}{}}
\rput(11,3){\rnode{Ax2}{}}
\rput(11,-37){\rnode{Bx2}{}}
              \rput(9,-26){\rnode{202}{$c_i\!\cdot\!\widetilde{\sigma}$}}
              \rput(0,-26){$\cdots$}
              \rput(-9,-26){\rnode{212}{$c_i\!\cdot\!\sigma$}}
              \rput(0,-16){\rnode{32}{$c_i\!\cdot\!\tau$}}      
              \rput(9,-6){\rnode{102}{$c_i\!\cdot\!\widetilde{\textbf{u}}$}}
              \rput(0,-6){$\cdots$}
              \rput(-9,-6){\rnode{112}{$c_i\!\cdot\!\textbf{u}$}}
              \end{tiny}}	
\rput(110,22){
\begin{tiny}
\rput(-7,20){$ \displaystyle \textbf{I}_{(i-1)}^{(n+1)}$\pscurve[linewidth=0.4pt,linestyle=solid]{->}(-10,-3)(-6,-7)(-4,-12)}  
\rput(0,9.6){\rnode{A3}{$c_{\textrm{i-1}}\!\cdot\!\omega_{n}$}}
\rput(0,-43.6){\rnode{B3}{}}
\rput(0,-46){$c_{i-1}$}
\rput(20,-17){\rnode{aa3}{}}
\rput(-20,-17){\rnode{bb3}{}}
\rput(-11,3){\rnode{A13}{}}
\rput(-11,-37){\rnode{B13}{}}
\rput(11,3){\rnode{Ax3}{}}
\rput(11,-37){\rnode{Bx3}{}}
              \rput(9,-26){\rnode{203}{$c_{\textrm{i-1}}\!\cdot\!\widetilde{\sigma}$}}
              \rput(0,-26){$\cdots$}
              \rput(-9,-26){\rnode{213}{$c_{\textrm{i-1}}\!\cdot\!\sigma$}}
              \rput(0,-16){\rnode{33}{$c_{\textrm{i-1}}\!\cdot\!\tau$}}      
              \rput(9,-6){\rnode{103}{$c_{\textrm{i-1}}\!\cdot\!\widetilde{\textbf{u}}$}}
              \rput(0,-6){$\cdots$}
              \rput(-9,-6){\rnode{113}{$c_{\textrm{i-1}}\!\cdot\!\textbf{u}$}}
              \end{tiny}}									
							\psset{nodesep=1pt,linestyle=solid,linewidth=0.4pt}  
								\ncline{-}{202}{203}
								\ncline{-}{212}{213}
								\ncline{-}{102}{103}
								\ncline{-}{112}{113}	
\rput(185,35){\rput(-50,13){\psline[linestyle=dotted,linewidth=1pt]{->}(0,0)(35,6)\rput{8}(17,6){\begin{tiny}$s_{1}\cdot  \textbf{x}$\end{tiny}}
}
\begin{tiny}
\rput(-7,20){$ \displaystyle \textbf{I}_{(1)}^{(n+1)}$\pscurve[linewidth=0.4pt,linestyle=solid]{->}(-10,-3)(-6,-7)(-4,-12)} 
\rput(0,9.6){\rnode{AL}{$\omega_{n+1}$}}

\rput(0,-43.6){\rnode{BL}{}}
\rput(0,-45){$c_1$}
\rput(20,-17){\rnode{aaL}{}}
\rput(-20,-17){\rnode{bbL}{}}
\rput(-11,3){\rnode{A1L}{}}
\rput(-11,-37){\rnode{B1L}{}}
\rput(11,3){\rnode{AxL}{}}
\rput(11,-37){\rnode{BxL}{}}
              \rput(9,-26){\rnode{20L}{$c_1\!\cdot\!\widetilde{\sigma}$}}
              \rput(0,-26){$\cdots$}
              \rput(-9,-26){\rnode{21L}{$c_1\!\cdot\!\sigma$}}
              \rput(0,-16){\rnode{3L}{$c_1\!\cdot\!\tau$}}      
              \rput(9,-6){\rnode{10L}{$c_1\!\cdot\!\widetilde{\textbf{u}}$}}
              \rput(0,-6){$\cdots$}
              \rput(-9,-6){\rnode{11L}{$c_1\!\cdot\!\textbf{u}$}}
							
\rput(-35,-7){\rput(0,9.6){\rnode{ALn}{}}
\rput(0,-43.6){\rnode{BLn}{}}
\rput(9,-26){\rnode{20Ln}{}}
              \rput(-9,-26){\rnode{21Ln}{}}
              \rput(0,-16){\rnode{3Ln}{}}      
              \rput(9,-6){\rnode{10Ln}{}}
              \rput(-9,-6){\rnode{11Ln}{}}
}
              \end{tiny}}		
						\psset{nodesep=0.5pt,linestyle=dashed,linewidth=0.4pt}  
								\ncline{-}{AL}{ALn}
								\ncline{-}{BL}{BLn}
								\ncline{-}{20Ln}{20L}
								\ncline{-}{21Ln}{21L}
								\ncline{-}{3Ln}{3L}
								\ncline{-}{10Ln}{10L}
								\ncline{-}{11Ln}{11L}			
									
									}
                                  
                      \psset{nodesep=1.3pt}
               
                 \psset{nodesep=0pt,offset=0pt} 
                 \ncline{-}{nol}{nur}
             \ncline{-}{nor}{nul} 
              \psset{nodesep=1.2pt} 
								\ncline{-}{A3}{A2}
								\ncline{-}{32}{33}
								  \ncline{-}{B2}{B3}
              \ncline{-}{3}{20}
             \ncline{-}{20}{3}
             \ncline{-}{3}{21}
             \ncline{-}{21}{3}
             \ncline{-}{3}{10}
             \ncline{-}{10}{3}
             \ncline{-}{3}{11}
             \ncline{-}{11}{3}
              \ncline{-}{A1}{A}
              \ncline{-}{B1}{B}
              \ncline{-}{B}{B1}
              \ncline{-}{A}{Ax}
              \ncline{-}{Ax}{A}
              \ncline{-}{B}{Bx}
              \ncarc[linewidth=0.5pt,arcangle=-15,linestyle=dashed]{-}{sn}{20}
              \ncarc[linewidth=0.5pt,arcangle=15,linestyle=dashed]{-}{ln}{10}
              \ncarc[arcangle=25,linestyle=dashed]{-}{A}{10}
              \ncarc[arcangle=-25,linestyle=dashed]{-}{A}{11}
              
              \ncarc[arcangle=-25,linestyle=dashed]{-}{B}{20}
              \ncarc[arcangle=25,linestyle=dashed]{-}{B}{21}
              
              \ncarc[arcangle=-25,linestyle=dashed]{-}{aa}{Ax}
              \ncarc[arcangle=25,linestyle=dashed]{-}{aa}{Bx}
              \ncarc[arcangle=-25,linestyle=dashed]{-}{bb}{B1}
              \ncarc[arcangle=25,linestyle=dashed]{-}{bb}{A1}
             \ncline{-}{202}{32}
             \ncline{-}{32}{212}
             \ncline{-}{212}{32}
             \ncline{-}{32}{102}
             \ncline{-}{102}{32}
             \ncline{-}{32}{112}
             \ncline{-}{112}{32}
              \ncline{-}{A2}{A12}
              \ncline{-}{A12}{A2}
              \ncline{-}{B12}{B2}
              \ncline{-}{B2}{B12}
              \ncline{-}{A2}{Ax2}
              \ncline{-}{Ax2}{A2}
              \ncline{-}{Bx2}{B2}
              \ncline{-}{B2}{Bx2}
              \ncarc[linewidth=0.5pt,arcangle=-15,linestyle=dashed]{-}{sn}{20}
              \ncarc[linewidth=0.5pt,arcangle=15,linestyle=dashed]{-}{ln}{10}
              \ncarc[arcangle=25,linestyle=dashed]{-}{A2}{102}
              \ncarc[arcangle=-25,linestyle=dashed]{-}{A2}{112}
              
              \ncarc[arcangle=-25,linestyle=dashed]{-}{B2}{202}
              \ncarc[arcangle=25,linestyle=dashed]{-}{B2}{212}
              
              \ncarc[arcangle=-25,linestyle=dashed]{-}{aa2}{Ax2}
              \ncarc[arcangle=25,linestyle=dashed]{-}{aa2}{Bx2}
              \ncarc[arcangle=-25,linestyle=dashed]{-}{bb2}{B12}
              \ncarc[arcangle=25,linestyle=dashed]{-}{bb2}{A12}							
             \ncline{-}{203}{33}
             \ncline{-}{33}{213}
             \ncline{-}{213}{33}
             \ncline{-}{33}{103}
             \ncline{-}{103}{33}
             \ncline{-}{33}{113}
             \ncline{-}{113}{33}
              \ncline{-}{A3}{A13}
              \ncline{-}{A13}{A3}
              \ncline{-}{B13}{B3}
              \ncline{-}{B3}{B13}
              \ncline{-}{A3}{Ax3}
              \ncline{-}{Ax3}{A3}
              \ncline{-}{Bx3}{B3}
              \ncline{-}{B3}{Bx3}
              \ncarc[linewidth=0.5pt,arcangle=-15,linestyle=dashed]{-}{sn}{20}
              \ncarc[linewidth=0.5pt,arcangle=15,linestyle=dashed]{-}{ln}{10}
              \ncarc[arcangle=25,linestyle=dashed]{-}{A3}{103}
              \ncarc[arcangle=-25,linestyle=dashed]{-}{A3}{113}
              
              \ncarc[arcangle=-25,linestyle=dashed]{-}{B3}{203}
              \ncarc[arcangle=25,linestyle=dashed]{-}{B3}{213}
              
              \ncarc[arcangle=-25,linestyle=dashed]{-}{aa3}{Ax3}
              \ncarc[arcangle=25,linestyle=dashed]{-}{aa3}{Bx3}
              \ncarc[arcangle=-25,linestyle=dashed]{-}{bb3}{B13}
              \ncarc[arcangle=25,linestyle=dashed]{-}{bb3}{A13}											
             \ncline{-}{20L}{3L}
             \ncline{-}{3L}{21L}
             \ncline{-}{21L}{3L}
             \ncline{-}{3L}{10L}
             \ncline{-}{10L}{3L}
             \ncline{-}{3L}{11L}
             \ncline{-}{11L}{3L}
              \ncline{-}{AL}{A1L}
              \ncline{-}{A1L}{AL}
              \ncline{-}{B1L}{BL}
              \ncline{-}{BL}{B1L}
              \ncline{-}{AL}{AxL}
              \ncline{-}{AxL}{AL}
              \ncline{-}{BxL}{BL}
              \ncline{-}{BL}{BxL}
              \ncarc[linewidth=0.5pt,arcangle=-15,linestyle=dashed]{-}{sn}{20}
              \ncarc[linewidth=0.5pt,arcangle=15,linestyle=dashed]{-}{ln}{10}
              \ncarc[arcangle=25,linestyle=dashed]{-}{AL}{10L}
              \ncarc[arcangle=-25,linestyle=dashed]{-}{AL}{11L}
              
              \ncarc[arcangle=-25,linestyle=dashed]{-}{BL}{20L}
              \ncarc[arcangle=25,linestyle=dashed]{-}{BL}{21L}
              
              \ncarc[arcangle=-25,linestyle=dashed]{-}{aaL}{AxL}
              \ncarc[arcangle=25,linestyle=dashed]{-}{aaL}{BxL}
              \ncarc[arcangle=-25,linestyle=dashed]{-}{bbL}{B1L}
              \ncarc[arcangle=25,linestyle=dashed]{-}{bbL}{A1L}										              
\end{pspicture}
\label{sym_n_as_sym_n+1}
\end{numeris}
 \textrm{  }\\[5.4cm]

Section ~\ref{proof_main_Theorem} is devoted to the proof of the following main theorem,which states that every first-degree critical permutation in $S_n$ gives rise to first-degree critical permutations in all $n+1$ left cosets, $\textbf{I}_{(i)}^{(n+1)}$, and right cosets, $\textbf{I}_{(n+1)}^{(i)}$, of $S_{n+1}$. Moreover, both the inverse and the $\omega_{n+1}$-conjugate of a critical permutation $\tau$ are also critical, due to the analogous Bruhat-order structure in $\Lambda_{(\tau)}$ (here $\omega_{n+1}$ denotes the longest element in $S_{n+1}$) (for more details, see Remark~\ref{isomorph}).

\begin{numeris}\normalfont{\textbf{Theorem.}} 
\textit{Let  $\tau \in S_{n}$ be considered  as  a permutation in  $S_{n+1}$ and let   $c_i$,  $d_i$ be  defined as  in \ref{sym_n_as_sym_n+1}. 
  If  $\tau $ is  critical in the first degree in $S_{n}$,  then}
	 
	 \begin{center}
	 \textit{$c_{i}\cdot  \tau$, \  \  \  \  \    $\tau \cdot   d_{i} $,   \  \  \  \  \   
		$\omega_{n+1}\cdot  \left(c_{i}\cdot  \tau\right)\cdot  \omega_{n+1}$,  \  \  \  \  \  $\omega_{n+1}\cdot \left(\tau \cdot   d_{i}\right)\cdot  \omega_{n+1}$, }
\end{center}
		
\textit{as well as their inverses,  are first-degree critical permutations  in $S_{n+1}$ for all $i\in \left[1,n+1\right]$.
}\label{Main_Theorem}
\end{numeris}

\begin{numeris}\normalfont{\textbf{Remark.}} Theorem~\ref{Main_Theorem}  does not generalize to higher degrees. For a critical permutation $\tau \in S_n$ of degree $\textbf{\textit{i}}$, the corresponding  permutations in the \textit{blocks} $\textbf{I}_{(i)}^{(n+1)}$ (or $\textbf{I}_{(n+1)}^{(i)}$) are not necessarily critical of degree $\textit{\textbf{i}}$. For example $\tau = (453126)\in S_5\subseteq S_6$ is critical of degree \textit{\textbf{2}}, but $c_4\cdot  \tau  = (563124)\in \textbf{I}_{(4)}^{(6)}$ is a critical   of degree \textit{\textbf{1}} in $S_6$.
 Nevertheless, for every critical permutation $\tau \in S_n$, the corresponding permutations $c_i \tau$ and $\tau d_i$ in the blocks remain critical, though their degree may differ. A full proof of this fact requires additional technical results and will be presented in future work.
Moreover, not every first-degree critical permutation in $S_{n+1}$ can be obtained from $S_n$ using the construction described in Theorem~\ref{Main_Theorem} (see Section~\ref{sec:OtherCriticalPermutations} for details).
\label{not_generalizable}
\end{numeris}

From this point onward,  let $\mathcal{C}_{n}$ denote the set of critical permutations in $S_4$, or to articulate this in terms of Verma multiplicities
\begin{center}
$\mathcal{C}_n := \left\{q\in S_n \mid q \textrm{ is a critical permutation} \right\} = \left\{q\in S_n \mid \left[M(\textrm{id}_n):L(q)\right]\geq 2\right\}$
\end{center}
We have already established in ~\ref{critical_example}, that $\mathcal{C}_1 = \mathcal{C}_2= \mathcal{C}_3 = \emptyset$, and  $\mathcal{C}_4 = \left\{ (3412) , \   (4231)\right\}$. Moreover,   $\textbf{x}_1 := (3412) $ and  $\textbf{x}_2 := (4231)$  are first-degree critical permutations.  
We proceed to determine all critical permutations in $S_{n+1}$ that arise from $S_n$ for $n=4,5,6$, with particular attention to those permutations that cannot be derived from  $\mathcal{C}_n$ using  the Theorem~\ref{Main_Theorem}.

\begin{numeris}\normalfont{\textbf{Permutations in $\mathcal{C}_{5}$.}} 
The diagram below symbolically   displays   permutations of $S_5$, organized based on their lengths (due to the impossibility of visualizing all neighboring lines without significant loss of clarity, all of them are disregarded). The permutations that do not have \textit{eye-catching} notation are non-critical permutations. The 32 permutation in $\mathcal{C}_{5}$  are depicted  using boxes, ellipses, and trapezoids. 
The left-right symmetry reflects the $\omega_5$-conjugacy-relation via the function $\Phi_5: S_5 \to S_5$ given by $\Phi_5(\sigma) =\omega_5\cdot  \sigma \cdot  \omega_5 $ between the depicted elements on opposite sides. For example  $l(34125)=l(\Phi_5(34125)) =5$, thus  $(34125)$ and $\Phi_5(34125) = (41523)$ are at 5 on the length scale,  and on opposite sides of the center of the diagram.
 For  $\tau$ on the dashed line we have $\tau = \Phi_5(\tau)$ (note that this does not apply to the entire diagram, since there are several permutations on the same length scale that are equal to their $\omega_5$-conjugates). 
The significance of the  permutations in ellipses, white, black and gray boxes is as follows:
 \begin{subnum}\normalfont{\textbf{Permutations in} $\mathcal{C}_{5}$ \textbf{by applying Theorem} ~\ref{Main_Theorem} \textbf{on} $\mathcal{C}_{4}$}. The permutations in  white boxes  originate from  $\textbf{x}_1 = (3412) $ and  in  ellipses from $\textbf{x}_2 = (4231)$,  by applying Theorem ~\ref{Main_Theorem}. For example \begin{scriptsize}$\underbrace{\framebox{(45132)}}_{c_2\cdot  \textbf{x}_1}$\end{scriptsize} indicates  $(45132) = \underbrace{s_2\cdot  s_3\cdot  s_4}_{c_2}\cdot  \underbrace{(34125)}_{\textbf{x}_1}$ and the meaning of \begin{scriptsize}$\underbrace{\ovalnode{B}{(42531)}}_{\textbf{x}_2\cdot  d_3}$\end{scriptsize} is $(42531) = \underbrace{(4231)}_{\textbf{x}_2}\cdot  \underbrace{s_4\cdot  s_3}_{d_3}$. 
The permutations in  gray boxes can be "\textit{derived}" from both $\textbf{x}_{1}$ and $\textbf{x}_{2}$ as indicated below the box. 
Thus 26 out of 32  elements in $\mathcal{C}_{5}$ are linked to $\textbf{x}_1$  and $\textbf{x}_2$.
\label{crit_in_S_5_from_S_4}  
\end{subnum}
\begin{subnum}\normalfont{\textbf{Permutations in} $\mathcal{C}_{5}$ \textbf{that do not come from} $\mathcal{C}_{4}$} (see Section~\ref{sec:OtherCriticalPermutations}). 
\begin{itemize}
	\item[$\bullet$]  $ \eta_n := \begin{scriptsize}\left(
\begin{array}{ccccccccccccccc}
	1 &2& \cdots & a-t & \cdots &a -2 &  a -1  & a & a+1 & a+2 & \cdots & a+ t & \cdots & n \\
	n-1 & n-3 & \cdots & 2\cdot t & \cdots &4 &  2 & n & 1 & 3 & \cdots & 2\cdot t-1 & \cdots & n-2
\end{array}
\right)\end{scriptsize}$ are first degree critical permutation  for all  $n=2a-1$ with $a\geq 3$. Thus $\eta_5=(42513)$, and its $\omega_5$-conjugate $(35142)$,  are   in trapezoids  at level 6 on the length  scale.  
	\item[$\bullet$] 	
	Permutations of the form $\tau = (n\ \textbf{a}\ \ldots\ \textbf{b}\ 1)$ with $\textbf{a}\neq n-1$ and $\textbf{b}\neq 2$ have at least five small neighbors for all $n\geq 5$. In $S_5$ there are exactly three such (self-inverse) cases, $(5\ \textbf{a}\ x\ \textbf{b}\ 1)$ with $\textbf{a}\neq 4$, $\textbf{b}\neq 2$, shown as \textit{reversed} trapezoids at levels 7 and 8.
	
\begin{center}
$\tau_{1}=(52341)$, \  \  \   $\tau_{2}=(52431)$, \  \  \   $\tau_{3}=(53241)$
\end{center}
\item[$\bullet$]   We will see in ~\ref{second_degree}  that,  for all $n= 2a-1$ with $a\geq 3$, the permutation of the form
$p_n =\begin{scriptsize}  \left(
\begin{array}{ccccccccc}
	1 & 2 & \cdots & a-1 & a & a+1& a+2 & \cdots & n\\
	a+1 & a+2 & \cdots & n & a & 1& 2 & \cdots & a-1
\end{array}
\right) \end{scriptsize}$
 is a critical of degree \textbf{\textit{2}}. For $a=3$ the second degree critical  permutation  $p_5= (45312)$ is located in the black box.
\end{itemize}
\label{crit_in_S_5_not_from_S_4}  
\end{subnum}

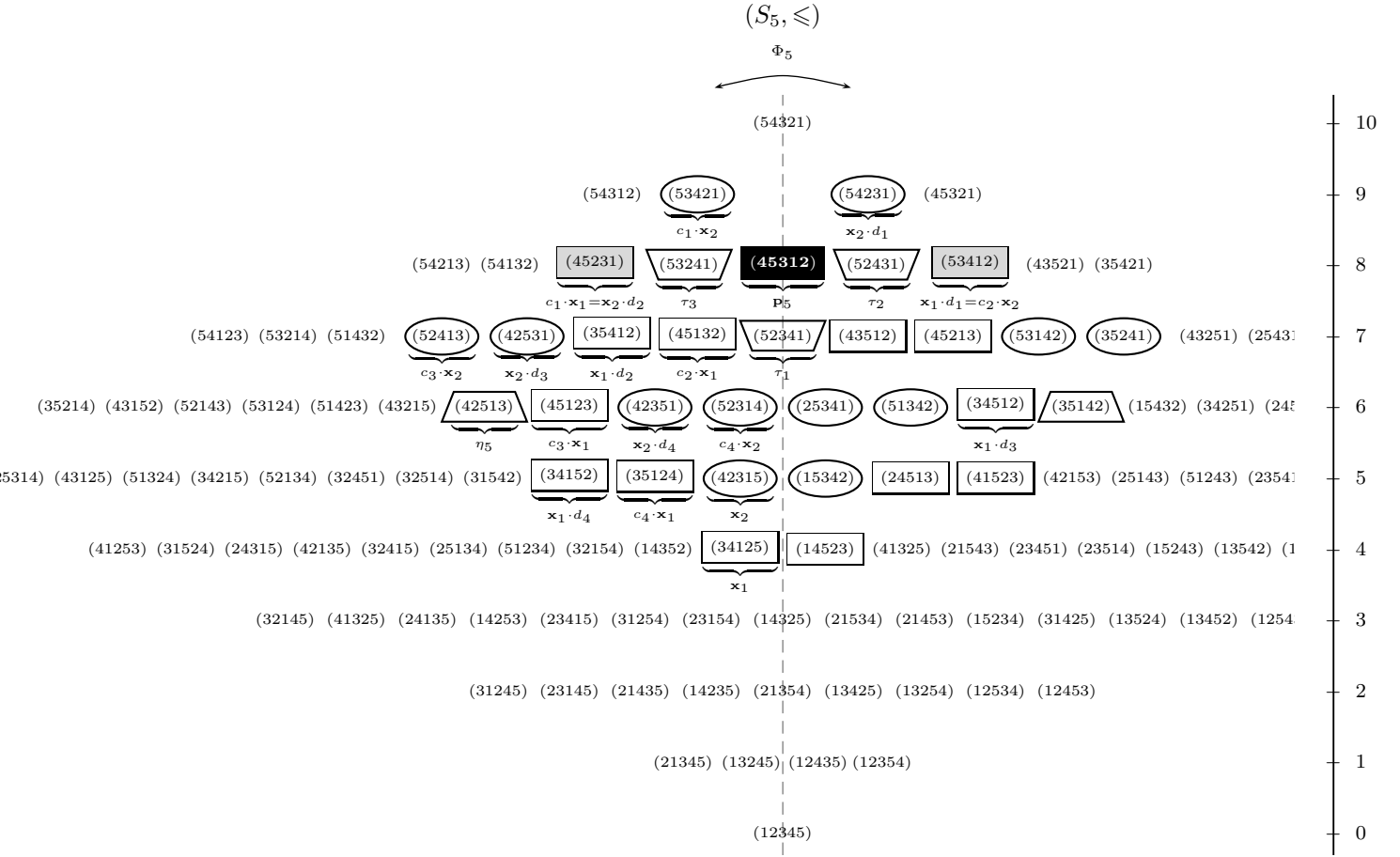
\begin{tiny}
\psset{xunit=2.4mm,yunit=1mm,runit=1mm}
\begin{pspicture}(0,0)(-35,113)
\rput(0,110){$\Phi_5$}
\rput(0,105){\pscurve[linewidth=0.5pt,linestyle=solid]{<->}(-4,0)(0,1.7)(4,0)}
\rput(0,115){\begin{footnotesize}$\left(S_5,\ma \right)$\end{footnotesize}}							
		\rput(0,-3){\rnode{obensig}{}} 
   \rput(0,104){\rnode{unten}{}}   					
\ncline[linewidth=0.5pt,linestyle=dashed,linecolor=gray]{obensig}{unten}
              \rput(0,0){\rnode{12345}{$(12345)$}}
							\psset{xunit=2mm,yunit=1mm,runit=1mm}
              \rput(-7,10){\rnode{21345}{$(21345)$}}
              \rput(-2.2,10){\rnode{13245}{$(13245)$}}  
              \rput(2.5,10){\rnode{12435}{$(12435)$}} 
              \rput(7,10){\rnode{12354}{$(12354)$}} 
              \rput(-20,20){\rnode{31245}{$(31245)$}} 
              \rput(-15,20){\rnode{23145}{$(23145)$}} 
              \rput(-10,20){\rnode{21435}{$(21435)$}} 
              \rput(-5,20){\rnode{14235}{$(14235)$}} 
              \rput(0,20){\rnode{21354}{$(21354)$}} 
              \rput(5,20){\rnode{13425}{$(13425)$}} 
              \rput(10,20){\rnode{13254}{$(13254)$}} 
              \rput(15,20){\rnode{12534}{$(12534)$}} 
              \rput(20,20){\rnode{12453}{$(12453)$}}
              \rput(-35,30){\rnode{32145}{$(32145)$}}  
              \rput(-30,30){\rnode{41235}{$(41325)$}} 
              \rput(-25,30){\rnode{24135}{$(24135)$}} 
              \rput(-20,30){\rnode{14253}{$(14253)$}} 
              \rput(-15,30){\rnode{23415}{$(23415)$}} 
              \rput(-10,30){\rnode{31254}{$(31254)$}} 
              \rput(-5,30){\rnode{23154}{$(23154)$}} 
              \rput(0,30){\rnode{14325}{$(14325)$}} 
              \rput(5,30){\rnode{21534}{$(21534)$}} 
              \rput(10,30){\rnode{21453}{$(21453)$}} 
              \rput(15,30){\rnode{15234}{$(15234)$}} 
              \rput(20,30){\rnode{31425}{$(31425)$}} 
              \rput(25,30){\rnode{13524}{$(13524)$}} 
              \rput(30,30){\rnode{13452}{$(13452)$}} 
              \rput(35,30){\rnode{12543}{$(12543)$}}
							\psset{xunit=2.4mm,yunit=1mm,runit=1mm}
              \rput(-39,40){\rnode{41253}{$(41253)$}}
              \rput(-35,40){\rnode{31524}{$(31524)$}} 
              \rput(-2.5,38){\rnode{34125}{$\underbrace{\frame{\colorbox{white}{(34125)}}}_{\textbf{x}_1}$}} 
              \rput(-27,40){\rnode{42135}{$(42135)$}} 
              \rput(-23,40){\rnode{32415}{$(32415)$}} 
              \rput(-19,40){\rnode{25134}{$(25134)$}} 
              \rput(-15,40){\rnode{51234}{$(51234)$}} 
              \rput(-11,40){\rnode{32154}{$(32154)$}} 
              \rput(-31,40){\rnode{24315}{$(24315)$}} 
              \rput(-7,40){\rnode{14352}{$(14352)$}} 
              \rput(7,40){\rnode{41325}{$(41325)$}} 
              \rput(31,40){\rnode{15324}{$(15324)$}} 
              \rput(11,40){\rnode{21543}{$(21543)$}} 
              \rput(15,40){\rnode{23451}{$(23451)$}} 
              \rput(19,40){\rnode{23514}{$(23514)$}} 
              \rput(23,40){\rnode{15243}{$(15243)$}} 
              \rput(27,40){\rnode{13542}{$(13542)$}} 
              \rput(2.5,40){\rnode{14523}{$\frame{\colorbox{white}{(14523)}}$}} 
              \rput(35,40){\rnode{24153}{$(24153)$}} 
              \rput(39,40){\rnode{31452}{$(31452)$}}  
              \rput(-41,50){\rnode{43125}{$(43125)$}}    
							  \rput(-2.5,48){\psellipse(0,2)(2.2,2.7)\rput(0,-1.34){$\underbrace{\textcolor{white}{----}}_{\textbf{x}_2}$} \rput(0,2){$(42315)$}}
              \rput(-37,50){\rnode{51324}{$(51324)$}}
              \rput(-33,50){\rnode{34215}{$(34215)$}}
              \rput(-29,50){\rnode{52134}{$(52134)$}}
              \rput(7.5,50){\rnode{24513}{$\frame{\colorbox{white}{(24513)}}$}}
              \rput(-21,50){\rnode{32514}{$(32514)$}}
              \rput(12.5,50){\rnode{41523}{$\frame{\colorbox{white}{(41523)}}$}}
              \rput(-17,50){\rnode{31542}{$(31542)$}}
              \rput(-25,50){\rnode{32451}{$(32451)$}}
              \rput(-45,50){\rnode{25314}{$(25314)$}}
              \rput(45,50){\rnode{41352}{$(41352)$}}
              \rput(25,50){\rnode{51243}{$(51243)$}}
              \rput(17,50){\rnode{42153}{$(42153)$}}
              \rput(-12.5,48){\rnode{34152}{$\underbrace{\frame{\colorbox{white}{(34152)}}}_{\textbf{x}_1\cdot  d_4}$}}
              \rput(21,50){\rnode{25143}{$(25143)$}}
              \rput(-7.5,48){\rnode{35124}{$\underbrace{\frame{\colorbox{white}{(35124)}}}_{c_4\cdot  \textbf{x}_1}$}}
              \rput(29,50){\rnode{23541}{$(23541)$}}
              \rput(33,50){\rnode{15423}{$(15423)$}}
              \rput(37,50){\rnode{24351}{$(24351)$}}
							\rput(2.5,50){\psellipse(0,0)(2.2,2.7)\rput(0,0){$(15342)$}}
              \rput(41,50){\rnode{14532}{$(14532)$}} 
              \rput(-42,60){\rnode{35214}{$(35214)$}}
							 \rput(-17.5,58){\pspolygon(-2.5,0)(2.5,0)(1.8,4.2)(-1.8,4.2)\rput(0,-1.34){$\underbrace{\textcolor{white}{----}}_{\eta_5}$} \rput(0,2){$(42513)$}}
              \rput(-34,60){\rnode{52143}{$(52143)$}}
              \rput(-30,60){\rnode{53124}{$(53124)$}}
              \rput(-26,60){\rnode{51423}{$(51423)$}}
              \rput(-38,60){\rnode{43152}{$(43152)$}}
              \rput(-22,60){\rnode{43215}{$(43215)$}}
              \rput(-12.5,58){\rnode{45123}{$\underbrace{\frame{\colorbox{white}{(45123)}}}_{c_3\cdot \textbf{x}_1}$}}
							\rput(7.5,60){\psellipse(0,0)(2.2,2.7)\rput(0,0){$(51342)$}}
							\rput(2.5,60){\psellipse(0,0)(2.2,2.7)\rput(0,0){$(25341)$}}
							 \rput(-2.5,58){\psellipse(0,2)(2.2,2.7)\rput(0,-1.34){$\underbrace{\textcolor{white}{----}}_{c_4\cdot \textbf{x}_2}$} \rput(0,2){$(52314)$}}
							\rput(-7.5,58){\psellipse(0,2)(2.2,2.7)\rput(0,-1.34){$\underbrace{\textcolor{white}{----}}_{\textbf{x}_2\cdot  d_4}$} \rput(0,2){$(42351)$}}
              \rput(12.5,58){\rnode{34512}{$\underbrace{\frame{\colorbox{white}{(34512)}}}_{\textbf{x}_1\cdot  d_3}$}}
              \rput(22,60){\rnode{15432}{$(15432)$}}
              \rput(38,60){\rnode{41532}{$(41532)$}}
              \rput(26,60){\rnode{34251}{$(34251)$}}
              \rput(30,60){\rnode{24531}{$(24531)$}}
              \rput(34,60){\rnode{32541}{$(32541)$}}
							\rput(17.5,60){\pspolygon(-2.5,-2)(2.5,-2)(1.8,2.2)(-1.8,2.2)\rput(0,0){$(35142)$}}
              \rput(42,60){\rnode{25413}{$(25413)$}}
              \rput(-33,70){\rnode{54123}{$(54123)$}}  
              \rput(-29,70){\rnode{53214}{$(53214)$}} 
							\rput(15,70){\psellipse(0,0)(2.2,2.7)\rput(0,0){$(53142)$}}
							\rput(20,70){\psellipse(0,0)(2.2,2.7)\rput(0,0){$(35241)$}}
              \rput(-25,70){\rnode{51432}{$(51432)$}}  
              \rput(10,70){\rnode{45213}{$\frame{\colorbox{white}{(45213)}}$}}  
              \rput(-5,68){\rnode{45132}{$\underbrace{\frame{\colorbox{white}{(45132)}}}_{c_2\cdot \textbf{x}_1}$}}  
								\rput(0,68){\pspolygon(-2.5,4.2)(2.5,4.2)(1.8,0)(-1.8,0)\rput(0,-1.34){$\underbrace{\textcolor{white}{----}}_{\tau_1}$} \rput(0,2){$(52341)$}}
              \rput(5,70){\rnode{43512}{$\frame{\colorbox{white}{(43512)}}$}}  
              \rput(-10,68){\rnode{35412}{$\underbrace{\frame{\colorbox{white}{(35412)}}}_{\textbf{x}_1\cdot  d_2}$}}  
              \rput(25,70){\rnode{43251}{$(43251)$}}  
						  \rput(-20,68){\psellipse(0,2)(2.2,2.7)\rput(0,-1.34){$\underbrace{\textcolor{white}{----}}_{ c_3 \cdot  \textbf{x}_2}$} \rput(0,2){$(52413)$}} 
							\rput(-15,68){\psellipse(0,2)(2.2,2.7)\rput(0,-1.34){$\underbrace{\textcolor{white}{----}}_{\textbf{x}_2\cdot  d_3}$} \rput(0,2){$(42531)$}} 
              \rput(29,70){\rnode{25431}{$(25431)$}}  
              \rput(33,70){\rnode{34521}{$(34521)$}} 
              \rput(-20,80){\rnode{54213}{$(54213)$}} 
              \rput(-16,80){\rnode{54132}{$(54132)$}} 
              \rput(11,78){\rnode{53412}{$\underbrace{\frame{\colorbox{gray!30}{(53412)}}}_{\textbf{x}_1\cdot  d_1=c_2\cdot \textbf{x}_2}$}} 
							\rput(-5.5,78){\pspolygon(-2.5,4.2)(2.5,4.2)(1.8,0)(-1.8,0)\rput(0,-1.34){$\underbrace{\textcolor{white}{----}}_{\tau_3}$} \rput(0,2){$(53241)$}}
              \rput(0,78){\rnode{45312}{$\underbrace{\frame{\colorbox{black}{\color{white}(\textbf{45312})}}}_{\textbf{p}_5}$}} 
							\rput(5.5,78){\pspolygon(-2.5,4.2)(2.5,4.2)(1.8,0)(-1.8,0)\rput(0,-1.34){$\underbrace{\textcolor{white}{----}}_{\tau_2}$} \rput(0,2){$(52431)$}}
              \rput(-11,78){\rnode{45231}{$\underbrace{\frame{\colorbox{gray!30}{(45231)}}}_{c_1\cdot  \textbf{x}_1=\textbf{x}_2\cdot  d_2}$}} 
              \rput(16,80){\rnode{43521}{$(43521)$}} 
              \rput(20,80){\rnode{35421}{$(35421)$}}
              \rput(-10,90){\rnode{54312}{$(54312)$}}
							\rput(5,88){\psellipse(0,2)(2.2,2.7)\rput(0,-1.34){$\underbrace{\textcolor{white}{----}}_{\textbf{x}_2\cdot  d_1}$} \rput(0,2){$(54231)$}} 
						\rput(-5,88){\psellipse(0,2)(2.2,2.7)\rput(0,-1.34){$\underbrace{\textcolor{white}{----}}_{c_1\cdot  \textbf{x}_{2}}$} \rput(0,2){$(53421)$}} 
              \rput(10,90){\rnode{45321}{$(45321)$}}
              \rput(0,100){\rnode{54321}{$(54321)$}}
							
							\rput(30,110){\psframe*[linecolor=white](0,0)(20,-120)} 
							\rput(-7,0){\psline[linewidth=0.7pt]{-}(39.3,-3)(39.3,104)
							\rput[l](38.8,100){\begin{scriptsize}$- \ \ 10 $\end{scriptsize}}
							\rput[l](38.8,90){\begin{scriptsize}$- \ \ 9$\end{scriptsize}}
							\rput[l](38.8,80){\begin{scriptsize}$- \ \ 8$\end{scriptsize}}
							\rput[l](38.8,70){\begin{scriptsize}$- \ \ 7$\end{scriptsize}}
							\rput[l](38.8,60){\begin{scriptsize}$- \ \ 6$\end{scriptsize}}
							\rput[l](38.8,50){\begin{scriptsize}$- \ \ 5$\end{scriptsize}}
							\rput[l](38.8,40){\begin{scriptsize}$- \ \ 4$\end{scriptsize}}
							\rput[l](38.8,30){\begin{scriptsize}$- \ \ 3 $\end{scriptsize}}
							\rput[l](38.8,20){\begin{scriptsize}$- \ \ 2 $\end{scriptsize}}
							\rput[l](38.8,10){\begin{scriptsize}$- \ \ 1 $\end{scriptsize}}
							\rput[l](38.8,0){\begin{scriptsize}$- \ \ 0 $\end{scriptsize}}
							}    
\end{pspicture}
\end{tiny}

\textrm{  }


\label{critical_26}
\end{numeris}


\begin{numeris}\normalfont{\textbf{Permutations in} $\mathcal{C}_{6}$.} 
 There are   354  permutations $\mathcal{C}_{6}$, of which  12 are  critical   of degree \textbf{\textit{2}} (all those permutations are of the form $c_{i}\cdot p_5$,  or its inverse, or $\omega_6$-conjugated  for $i\in \left\{2,3,5,6\right\}$), \textit{one} is critical of degree \textbf{\textit{3}}, namely  $q_6=(564312)$, and remainder are of degree \textbf{\textit{1}}. For  
 \begin{center}
 $\ger{X} :=\left\{c_i\cdot q , \ q\cdot d_i \mid q\in \mathcal{C}_{5}, \ and \ i\in \left[1,6\right]\right\}$,
 \end{center}
	we have $\ger{Y}:= \ger{X}  \cup \left(\omega_6\cdot \ger{X} \cdot\omega_6\right) \subseteq \mathcal{C}_{6}$ with $\left|\mathcal{C}_{6}\backslash \ger{Y}\right|=11$ (see \cite{Elbos}). In other words, 11 permutations in $\mathcal{C}_6$ cannot be obtained from $\mathcal{C}_5$ via the construction described in the Theorem~\ref{Main_Theorem} (all of them are self-inverse). Except for $q_6$, the remaining 10 permutations are critical of degree \textbf{\textit{1}}. 
	
	It is easy to compute that there are 8 self-inverse permutations in  $S_{6}$ of the form 
 $ \left(6 \  \textbf{a} \ x \ y \  \textbf{b} \ 1\right)$ with $\textbf{a}\neq 5$ and  $\textbf{b}\neq 2$. All of them are in $\mathcal{C}_6$. 
 With the exception of $(645231)= c_1 \cdot  (53412)$, all others 7 permutations  do not arise from  $\mathcal{C}_5$. The remaining three permutations, $(361542)$, $(532614)$, and $(426153) \in \mathcal{C}_{6}$, are likewise not linked to $\mathcal{C}_5$ by Theorem~\ref{Main_Theorem}.
\label{Critical_in_S_6}
\end{numeris}

\begin{numeris}\normalfont{\textbf{Permutations in} $\mathcal{C}_{7}$.} Out of the 3\:488 critical permutations in $S_7$, there are 96 critical of degree \textbf{\textit{2}}, 12 are critical of degree \textbf{\textit{3}} and one is critical of degree \textbf{\textit{4}}.  
The  elements in  $\mathcal{C}_{7}$ exhibits similar patterns to the ones illustrated above for $\mathcal{C}_{5}$ and $\mathcal{C}_{6}$. There are 
only 20 permutations in $\mathcal{C}_{7}$ that are not linked to those in $\mathcal{C}_{6}$ (all of them are self-inverse). Among those,  15 are of the form 
 $ \left(7 \  \textbf{a} \ x \ y \ z \  \textbf{b} \ 1\right)$ with $\textbf{a}\neq 6$ and  $\textbf{b}\neq 2$.  The remaining five permutations are $(3716542)$, $(4271653)$, as well as they $\omega_7$-conjugates, and  $q_7=(6754312)$ that  is critical  in degree \textit{\textbf{4}}.
\label{Critical_in_S_7}
\end{numeris}

\section{First-Degree Critical Permutations}
\label{small_N_and_large_simples}
\begin{small}In this section, we provide the necessary details for the analysis of 
first-degree critical permutations.  We discuss structural properties of permutations 
in terms of the simple transpositions that appear in reduced decompositions of $\tau$, and we review some results concerning the  inversions that give rise to the \textit{small neighbors} 
of  $\tau$ with respect to the Bruhat order, which play a decisive role in determining its 
"criticality" (see, for example, \cite{Bjor},\cite{Fed}, \cite{Ten} for more details).\\
\end{small}

We begin by recalling the key notations and reviewing the context established in the previous sections: According to the Definition~\ref{Def_of_n_degree}, a permutation,  $\tau\in S_n$  is \textit{first-degree critical}   if and only if 
the number  of  the set of simple 
transpositions that are comparable with $\tau$, 
denoted by  
\begin{center}
$\textbf{L}_{\textbf{t}}(\tau) = \left\{s \in S_n \mid s \ma \tau, \ and \ l(s)=1 \right\}$.
\end{center}
 differs from the number of the set of  \textit{small neighbors} of $\tau$ with respect to Bruhat order (see ~\ref{Bruhat}), denoted by the set  $\textbf{N}_{\textbf{s}}(\tau)$. 
Using the fact that $\tau \cdot (a,b) = (c,d)\cdot\tau $  for transpositions $(a,b)$ and $(c,d)$ with  $(c,d) =(\tau(a),\tau(b))$ (see ~\ref{non_com}),  the set of small neighbors of $\tau$ is given by
\begin{center}
$
\begin{array}{ccl}
	\textbf{N}_{\textbf{s}}(\tau) &=&\left\{\sigma \in S_n \mid \sigma \sma \tau \right\} \\[1mm]
	&=& \left\{\sigma \in S_n \mid  l(\tau) = l(\sigma)+1, \ \textrm{and} \  \sigma = \tau\cdot (a,b), \ \textrm{where} \  (a,b)\in S_n \right\} \\[1mm]
	&=&  \left\{\sigma \in S_n \mid  l(\tau) = l(\sigma)+1, \ \textrm{and} \  \sigma = (a,b) \cdot \tau , \ \textrm{where} \ (a,b)\in S_n \right\}
\end{array}
$
\end{center}
Thus,  $\tau \in S_n$ is \textit{first} degree critical if and only if $\left|\textbf{N}_{\textbf{s}}(\tau)\right|\neq \left|\textbf{L}_{\textbf{t}}(\tau)\right|$.

\begin{numeris}\begin{center} \textbf{Reduced decompositions and free-fall inversions}\end{center}\end{numeris}
\begin{small}In \cite{Fed}, Incitti studied the comparability of permutations with respect to the Bruhat order. Building on these results, we derive the relationship between  $\textbf{N}_{\textbf{s}}(\tau)$ and a special class of inversions of $\tau$, as well as the connection between $\textbf{L}_{\textbf{t}}(\tau)$ and the simple transposition in a reduced decomposition of $\tau$.
\end{small}

\begin{subnum}\normalfont{\textbf{Reduced decomposition.}}  Because of the Coxeter relations, $s_i\cdot  s_j = s_j\cdot  s_i$ if $\left|i-j\right|\geq 2$, and $s_i\cdot  s_{i+1}\cdot  s_{i} = s_{i+1}\cdot  s_i\cdot  s_{i+1}$,  every  simple transposition $s_{i_j}$ that appears in one reduced decomposition of $\tau = s_{i_1} \cdot s_{i_2} \cdot \ldots \cdot s_{i_k}$, 
also appears in every reduced decomposition of $\tau$ (though not necessarily with the same multiplicity). In this paper we denote by $\Red(\tau)$ the set of all (pairwise different) simple transpositions that appear in a reduced decomposition of $\tau$: 
\begin{center}
$\Red(\tau) = \left\{ s_{i_1} , s_{i_2} ,\ldots , s_{i_k}\right\}$
\end{center}
Note that $\left|\Red(\tau)\right| \leq \left|\Inv(\tau)\right| = l(\tau)$ (see ~\ref{inversions}).  Inversions and reduced decompositions of permutations have been extensively studied (see \cite{Bjor}). In what follows, we recall these results and reformulate them in a way that aligns with our present approach.
\label{reduced_decom}
\end{subnum}

	\begin{subnum}\normalfont{\textbf{Lemma.}} \textit{Let  $\eta\in S_n$. }
	\begin{itemize}
		\item[(1)] \textit{$\textbf{L}_{\textbf{t}}(\eta)  =  \Red(\eta)$.} 
		\item[(2)] \textit{If   $s$ is a  a simple transposition and  $s\not\in \Red(\eta)$, then}	 
		\begin{center}
		\textit{$\left|\textbf{L}_{\textbf{t}}(s\cdot\eta)\right| = \left|\textbf{L}_{\textbf{t}}(\eta)\right| +1$ \  \  \  \  and  \  \  \  \  
						$\left|\textbf{N}_{\textbf{s}}(s\cdot\eta)\right| = \left|\textbf{N}_{\textbf{s}}(\eta)\right| +1$}
\end{center}
		 \item[(3)]  \textit{$\left|\textbf{N}_{\textbf{s}}(\eta)\right| \geq \left|\textbf{L}_{\textbf{t}}(\eta)\right| $.}
	\item[(4)]  \textit{$\eta$ is critical of degree one if and only if  $\left|\textbf{N}_{\textbf{s}}(\eta)\right|  > \left|\textbf{L}_{\textbf{t}}(\eta)\right| $.}

	\end{itemize}
	
	\textit{Proof.} (1) A simple transposition $s$ appears in   $\Red(\eta)$ if and only if $s \ma \eta$ with respect to the Bruhat  order (see \cite[Corollary 2.2.3]{Bjor}). Thus  $\textbf{L}_{\textbf{t}}(\eta)  \supseteq  \Red(\eta)$  and  $\textbf{L}_{\textbf{t}}(\eta) \subseteq  \Red(\eta)$. 
	
	(2)  If   $s\not\in \Red(\eta)$, then  
by \cite[Proposition 2.2.7]{Bjor},  for $x\in S_n$ with $x\ma s\cdot \eta$, we have $x=\eta$ or   $x=s\cdot \sigma$ for some  $\sigma\ma\eta$. 
 Since $l(s\cdot\eta)=l(\eta)+1$ and for each $\sigma$ with $\sigma\sma \eta $ we have $l(s\cdot\sigma)=l(\sigma)+1 = l(\eta)$. Thus, 
$\textbf{N}_{\textbf{s}}(c\cdot  \eta) = \left\{\eta\right\}\dot{\cup} \left\{s\cdot \sigma \mid \sigma \in \textbf{N}_{\textbf{s}}(\eta) \right\} $
  and therefore  $\left|\textbf{N}_{\textbf{s}}(s\cdot  \eta)\right| = \left|\textbf{N}_{\textbf{s}}(\eta)\right| +1$. 
Moreover, 	$ \Red(s\cdot\eta) =  \Red(\eta) \dot{\cup} \left\{s\right\}$, using (1) we have $\left|\textbf{L}_{\textbf{t}}(s\cdot\eta)\right| =\left|\Red(s\cdot\eta)\right| = \left| \Red(\eta) \dot{\cup} \left\{s\right\}\right| = \left|\textbf{L}_{\textbf{t}}(\eta)\right| +1$.

 (3) Let $m:=  \left|\textbf{L}_{\textbf{t}}(\eta)\right| = \left|\Red(\eta)\right|$. We  show $\left|\textbf{N}_{\textbf{s}}(\eta)\right| \geq m$  by induction on $m$: If $m=1$, then $\eta =s_i $, thus $\left|\textbf{N}_{\textbf{s}}(\eta)\right| = \left| \left\{\textrm{id}_n\right\}\right|  = 1 \geq 1$. 
 Let  $\eta =s\cdot  s_{i_1}\cdot  s_{i_2}\cdot  \ldots \cdot  s_{i_k}$ be  reduced  with $\sigma = s_{i_1}\cdot  s_{i_2}\cdot  \ldots \cdot  s_{i_k}$ and  $\left|\Red(\sigma)\right|=m$. By the assumption we have $\left|\textbf{N}_{\textbf{s}}(\sigma)\right| \geq m$. 
If $s\not\in \Red(\sigma)$, then by (2) we have $\left|\Red(\eta)\right| = \left|\Red(\sigma)\right| +1= m +1$ and  $\left|\textbf{N}_{\textbf{s}}(\eta)\right| = \left|\textbf{N}_{\textbf{s}}(s\cdot\sigma)\right| = \left|\textbf{N}_{\textbf{s}}(\sigma)\right| +1 \geq m+1$. And in case 
 $s\in \Red(\sigma)$, then $\left|\Red(\eta)\right| = \left|\Red(\sigma)\right| = m$, thus by assumption $\left|\textbf{N}_{\textbf{s}}(\eta)\right| \geq m$. 
	
	(4) Based on (3) we have $\left|\textbf{N}_{\textbf{s}}(\eta)\right| \neq \left|\textbf{L}_{\textbf{t}}(\eta)\right| $ if and only if $\left|\textbf{N}_{\textbf{s}}(\eta)\right| >  \left|\textbf{L}_{\textbf{t}}(\eta)\right| $.
	\hfill $\Box$ 
\label{simple_transpositions_of_tau}
\end{subnum}

A convenient characterization for adjacent permutations of $(S_n, \leqslant)$ is  introduce by  Incitti referring to the term  \textit{free-rise} (see \cite[Def. 2.3]{Fed}).   
Since  the number of small neighbors is crucial to determine critical permutations of degree one, we reformulate this concept using the term \textit{free-fall}. This allows us to determine the small neighbors of the permutation, providing a more convenient characterization.

\begin{subnum}\normalfont{\textbf{Definition.}} An inversion $(a,b)\in \Inv(\tau) := \left\{(a,b) \mid a<b \textrm{ and } \tau(a)>\tau(b)\right\}$ of $\tau$ is called  a \textit{\textit{free-fall}} inversion  of $\tau$ if there is \textit{no} $k\in \left[a,b\right]$ with $\tau(a)>\tau(k)>\tau(b)$.

  We denote by  $\Inv_{\textbf{f}}(\tau)$ the set of all free-fall inversions.  
\label{free_fall}
\end{subnum}

\begin{subnum}\normalfont{\textbf{Remark.}} 
By the properties of Coxeter relations, every permutation 
has  a reduced decomposition that can be grouped into \textit{blocks} via $\tau = \ger{b}_1 \cdot \ger{b}_2 \cdots \ger{b}_m$, where each \textit{block} $\ger{b}_i$ is a product of simple transpositions with consecutive indices, namely  $\Red(\ger{b}_i) = \left\{ s_j \mid j\in \left[l_i, k_i-1\right]\right\}$, and the indices of successive blocks strictly increase: In other words,  if $a<b$, then $t_a< t_b$ for all $t_a\in \left[l_a ,k_a-1\right]$ and $t_b\in \left[l_b, k_b-1\right]$. For example,
 $x =  s_{12}\cdot   s_5\cdot s_{9}\cdot  s_4\cdot s_6\cdot s_{10}  \cdot s_9 \cdot s_5$ is
\\[3mm]
$ x =   \underbrace{s_5\!\cdot\! s_4\!\cdot\! s_6\!\cdot\! s_5}_{\ger{b}_1}\!\cdot\! \underbrace{s_{9}\!\cdot\! s_{10}\!\cdot\! s_9}_{\ger{b}_2}\!\cdot\! \underbrace{s_{12}}_{\ger{b}_3}$   \    with  \    $\Red(\ger{b}_1)=\left\{s_{4}, s_{5},  s_{6} \right\}$,   \   $\Red(\ger{b}_2)=\left\{s_{9}, s_{10}\right\}$,    \  $\Red(\ger{b}_3)=\left\{s_{12}\right\}$.
\\[2mm]
Here $m=3$ with  $l_1=4$ and  $k_1= 7$, $l_2=9$ and  $k_2 = 11$, $l_3 = 12$ and $k_3=13$. 

Note that in the bottom line of $\tau$, only the numbers from  the intervals $[l_i, k_i]$ are permuted, while for every $t$ outside these intervals we have $\tau(t) = t$. For instance, in the example above 
\begin{center}
$x =  \begin{scriptsize}\Bigg(
\begin{array}{ccc}
	1 & 2 & 3 \\
	1&2&  3
\end{array}
\underbrace{\framebox{$
\begin{array}{cccc}
	4 & 5 & 6&7 \\
	6 & 7 & 4&5
\end{array}$
}}_{\ger{b}_1}
\begin{array}{c}
	8\\
	8 
\end{array}
\underbrace{\framebox{$
\begin{array}{ccc}
	9 & 10 & 11  \\
	11 & 10 & 9 
\end{array}$
}}_{\ger{b}_2}\underbrace{\framebox{$
\begin{array}{cc}
	12 & 13   \\
	13 & 12 
\end{array}$
}}_{\ger{b}_3}
\begin{array}{cccc}
	14 & \cdots &n&n+1\\
	14 & \cdots &n&n+1
\end{array} \Bigg)\end{scriptsize}$
\end{center}
It is easy to see, that if $a$ and $b$ lie in different intervals $[l_i, k_i]$, or if either of them lies outside all intervals, then $(a,b)$ is not an inversion of $\tau$, since $a<b$ implies $\tau(a)<\tau(b)$. Consequently, $(a,b)$ is an inversion of $\tau$ if and only if it is an inversion of $\ger{b}_i$ for some $i$. The same statement holds for the set of free-fall inversions.  Moreover, since for any $(a,b) \in \Inv_{\textbf{f}}(\tau)(\ger{b}_i)$ and $(c,d) \in \Inv_{\textbf{f}}(\tau)(\ger{b}_j)$ with $i \neq j$ we have $\{a,b\} \cap \{c,d\} = \emptyset$, the set of  free-fall inversions of $\tau$ is the disjoint union
$
\Inv_{\textbf{f}}(\tau) =  \Inv_{\textbf{f}}(\ger{b}_i) \cup \Inv_{\textbf{f}}(\ger{b}_2) \cup \cdots \cup \Inv_{\textbf{f}}(\ger{b}_m) 
$. For the permutation $x$ in the preceding example
\begin{center}
$\Inv_{\textbf{f}}(x) =  \underbrace{\left\{(4,6), (4,7), (5,6), (5,7)\right\}}_{\Inv_{\textbf{f}}(\ger{b}_1)} \cup \underbrace{\left\{(9,10), (10,11)\right\}}_{\Inv_{\textbf{f}}(\ger{b}_2)} \cup \underbrace{\left\{(12,13)\right\}}_{\Inv_{\textbf{f}}(\ger{b}_3)} .$
\end{center}
\label{Remark_about_blocks_general}
\end{subnum}

The following lemma is derived from \cite[Proposition 2.4.]{Fed}, reformulated in terms of the notion  \textit{free-fall} inversions, and   shows that
$\sigma \sma \tau$ if and only if $\sigma = \tau \cdot (a,b)$ for $(a,b)\in \Inv_{\textbf{f}}(\tau)$.

\begin{subnum}\normalfont{\textbf{Lemma.}}  \textit{If  $\tau\in (S_{n},\ma)$, then  the set of small neighbors of 
$\tau $ is precisely
}
\begin{center}
$\textbf{N}_{\textbf{s}}(\tau) = \left\{\tau\cdot (a,b) \mid (a,b)\in \Inv_{\textbf{f}}(\tau)\right\}$,
\end{center}
\textit{and consequently $\left|\textbf{N}_{\textbf{s}}(\tau)\right| = \left|\Inv_{\textbf{f}}(\tau)\right|$. }
\\

\textit{Proof.} Let $\sigma,\tau\in S_{n}$. According to  \cite[Proposition 2.4.]{Fed} we have  $\sigma\sma \tau $  if and only if $\sigma = \tau \cdot  (a,b)$ where $(a,b)$ is a free-rise of $\sigma$. The fact that $(a,b)$ is a free-rise of  $\sigma$, implies that $(a,b)$ is a free-fall of $\tau$. Thus $\sigma\in \textbf{N}_{\textbf{s}}(\tau) $ if and only if $\sigma = \tau \cdot (a,b)$ where $(a,b)\in \Inv_{\textbf{f}}(\tau).$
\hfill $\Box$ 
\label{neigh_lemma}
\end{subnum}

Building on the example from ~\ref{Remark_about_blocks_general}, we observe that the   $x$ has seven free-fall inversions, and therefore seven small neighbors. Combining Lemma~\ref{simple_transpositions_of_tau} with Lemma~\ref{neigh_lemma}, we obtain 
\begin{center}
$\left|\textbf{N}_{\textbf{s}}(x)\right| = 7$ \  \  \  \  and \  \  \  \  $\left|\textbf{L}_{\textbf{t}}(x)\right|=6$
\end{center}
 Hence, $x$ constitutes a first-degree critical permutation in $S_n$ for all $n \geq 13$.

\begin{subnum}\normalfont{\textbf{Remark.}} Observe that if $\tau(a) = j+1$ and $\tau(b) = j$ for some $a<b$, then in the  notation of $\tau$ the entries at positions $a$ and $b$ appear as $\tau = \begin{scriptsize} \left(
\begin{array}{ccccccc}
1 & \cdots  & a &  \cdots &b & \cdots & n\\
	\ast & \cdots & j+1 & \cdots & j& \cdots &\ast
\end{array}\right)
  \end{scriptsize}$. In this case, it is immediate that $(a,b)$ is a free-fall inversion of $\tau$, since no $k \in [a,b]$ can satisfy $\tau(a)>\tau(k)>\tau(b)$. Moreover, because   $\tau \cdot (a,b) = (\tau(a), \tau(b))\cdot \tau = (j+1,j)\cdot \tau = s_{j}\cdot \tau  $ and $\tau \cdot (a,b) \sma \tau$, it follows that 
	$s_j \cdot \tau \sma \tau$.
\label{free_fall_Remark}
\end{subnum}

\begin{numeris}\begin{center}\textbf{Bruhat Order on Cosets of} $S_{n+1}$ \textbf{over} $S_{n}$\end{center}\end{numeris}
\begin{small}Furthermore, we consider  $S_{n+1}$ as the union of  left  (or right) cosets of a subgroup $S_n$.
Following the notations in ~\ref{sym_n_as_sym_n+1}, we have $c_{i} \cdot  S_n = \textbf{I}_{(i)}^{(n+1)}$ and  $ S_n \cdot  d_{i} = \textbf{I}_{(n+1)}^{(i)}$, where   $c_{n+1} = d_{n+1} = \textrm{id}_{n+1}$ and $c_{i}= s_{i}\cdot  s_{i+1}\cdots  s_n$ as well as $d_{i}= \left(c_{i}\right)^{-1}$ for  $i\in \left[1,n\right]$. Recall that $\eta\in \textbf{I}_{(i)}^{(n+1)}$ if and only if $\eta(n+1)=i$, and  $\left|S_n\right|=\left|\textbf{I}_{(i)}^{(n+1)}\right| = \left|\textbf{I}_{(n+1)}^{(i)}\right|$.  We identify the elements of  $S_n$ with the elements of   $\textbf{I}_{(n+1)}^{(n+1)}$.  
\end{small}\\

The next Lemma shows that for each $\textbf{x}\in \textbf{I}_{(i-1)}^{(n+1)}$, the permutation $s_{i-1}\cdot  \textbf{x}$ is a uniquely determined small neighbor in $\textbf{I}_{(i)}^{(n+1)}$, which implies that in the diagram of $S_{n+1}$  each permutation in  $\textbf{I}_{(i-1)}^{(n+1)}$ is connected by a downward line to exactly one permutation in $\textbf{I}_{(i)}^{(n+1)}$. Moreover,   the functions  $S_n \xrightarrow{c_{i}\cdot  -} \textbf{I}_{(i)}^{(n+1)}$ are $S_n \xrightarrow{ -\cdot  d_{i}} \textbf{I}^{(i)}_{(n+1)}$,  preserve Bruhat-order, which yield that the shapes of diagrams of $(S_{n},\leqslant )$,  $(\textbf{I}_{(i)}^{(n+1)} , \leqslant )$ and $(\textbf{I}^{(i)}_{(n+1)} , \leqslant )$  have the same structural characteristic. Thus it shows that the Bruhat diagram of $S_{n+1}$ visualizes the neighboring-relationships both between blocks of the form $\textbf{I}_{(i)}^{(n+1)}$ (and $\textbf{I}_{(n+1)}^{(i)}$) as well as within each block, as illustrated in  ~\ref{sym_n_as_sym_n+1}.

\begin{subnum}\normalfont{\textbf{Lemma.}}  \textit{Let  $\sigma, \tau \in  S_n=\textbf{I}_{(n+1)}^{(n+1)}$, and  $c_i$,  $d_i$  as defined in ~\ref{sym_n_as_sym_n+1}, here  $i\in \left[1,n+1\right]$. }
\begin{itemize}
\item[(1)] \textit{ We have 
$ \left(c_{i}\cdot  \tau \right)\sma \left(c_{i-1}\cdot  \tau \right) $, \   and \  
	  $ \left(\tau \cdot  d_{i}\right)\sma \left(\tau  \cdot  d_{i-1}\right) $ for each  $i\in \left[2,n+1\right]$.}
	\item[(2)]  \textit{The following statements are equivalent}
\begin{center}
(i) $ \sigma \sma \tau$, \  \  \  \  \  \   \  \  \ 
	(ii)  $\left(c_{i}\cdot  \sigma\right) \sma \left(c_{i}\cdot  \tau\right)$, \  \  \  \  \  \  \  \  \ 
	iii)  $ \left(\sigma \cdot  d_{i}\right)\sma \left(\tau \cdot  d_{i}\right)$.
\end{center}
\end{itemize}

\textit{Proof.} (1) Consider $\textbf{x}:=c_{i-1}\cdot  \tau = \begin{scriptsize} \left(
\begin{array}{ccccc}
1 & \cdots  & l &  \cdots & n+1\\
	\ast & \cdots & i & \cdots & i-1 
\end{array}\right)
  \end{scriptsize} \in \textbf{I}_{(i-1)}^{(n+1)}$. We have  $s_{i-1} \cdot  \textbf{x} \sma \textbf{x}$ (see ~\ref{free_fall_Remark}).
	Since $s_{i-1}\cdot s_{i-1} = id_{n+1}$ and $c_i=s_{i}\cdot s_{i+1}\cdots s_n$, we have 
	$s_{i-1}\cdot  \textbf{x} = s_{i-1}\cdot  \left(c_{i-1}\cdot  \tau\right) = s_{i-1}\cdot  s_{i-1}\cdot s_{i}\cdot s_{i+1}\cdots s_n\cdot \tau  = c_{i}\cdot  \tau $. Thus 
 we have $ \left(c_{i}\cdot  \tau \right)\sma \left(c_{i-1}\cdot  \tau \right) $. 
	
	Recall that for all $\eta, \mu \in S_n$ with $\eta\Bruh \mu$ we also have $\eta^{-1}\Bruh \mu^{-1}$ (see \cite{Humph}). 
	Since $\tau\in \textbf{I}_{(n+1)}^{(n+1)}$ implies $\tau^{-1}\in \textbf{I}_{(n+1)}^{(n+1)}$, so as proved above  $ \left(c_{i}\cdot  \tau^{-1} \right)\sma \left(c_{i-1}\cdot  \tau^{-1} \right) $ for all $i\in \left[2,n+1\right]$. 
 According to the notations in  ~\ref{sym_n_as_sym_n+1} we have $\tau \cdot  d_{i} = \left(c_{i}\cdot  \tau^{-1} \right)^{-1} $ and $\tau \cdot  d_{i-1} = \left(c_{i-1}\cdot  \tau^{-1} \right)^{-1}$. Therefore 	$ \left(c_{i}\cdot  \tau^{-1} \right)\sma \left(c_{i-1}\cdot  \tau^{-1} \right) $  implies $  \left(\tau \cdot  d_{i}\right)\sma \left(\tau  \cdot  d_{i-1}\right)$.

	(2)  The part (1) implies $\eta\sma (c_{n}\cdot \eta)\sma (c_{n-1}\cdot  \eta) \sma \cdots \sma (c_{i+1}\cdot  \eta) \sma (c_{i}\cdot  \eta) \sma (c_{i-1}\cdot  \eta)$ for all $\eta\in S_n$. Thus   $l(c_{i}\cdot  \eta) = (c_{i+1}\cdot  \eta)+1$ for all $i\in \left[1,n\right]$ (see ~\ref{Bruhat}). By induction, we obtain $l(c_{i}\cdot  \eta) = l(\eta)+n+1-i$.
	
	"(i) $\Rightarrow$ (ii)" Since $\sigma, \tau \in S_n$ with $\sigma \sma \tau $, we have $\sigma = \tau \cdot  (a,b)$ and $l(\tau) = l(\sigma)+1$. The equation $\sigma = \tau \cdot  (a,b)$
	implies $c_i\cdot  \sigma = \left(c_i\cdot  \tau\right)\cdot  (a,b)$. Since $\sigma,\tau\in S_n$ we have   $l(c_{i}\cdot  \tau) = l(\tau)+n+1-i$, and  $l(c_{i}\cdot  \sigma) = l(\sigma)+n+1-i$. Thus   $l(\tau) = l(\sigma)+1$, imply that  $l(c_{i}\cdot  \tau)  = l(c_{i}\cdot  \sigma) +1 $. By the definition of Bruhat order, we have $\left(c_{i}\cdot  \sigma\right) \sma \left(c_{i}\cdot  \tau\right)$.
	
	"(ii) $\Rightarrow$ (i)" Since $\sigma,\tau\in S_n$ we have $l(c_i\cdot\sigma)=l(\sigma)+n+1-i$ and $l(c_i\cdot\tau)=l(\tau)+n+1-i$.   If $\left(c_{i}\cdot  \sigma\right) \sma \left(c_{i}\cdot  \tau\right)$, then $l\left(c_{i}\cdot  \tau\right) = l\left(c_{i}\cdot  \sigma\right) +1 $  and $c_{i}\cdot  \sigma = \left(c_{i}\cdot  \tau\right)\cdot  (a,b)$. Thus $l(\tau)+n+1-i= l\left(c_{i}\cdot  \tau\right) = l\left(c_{i}\cdot  \sigma\right) +1 = l(\sigma)+n+2-i$ implies  $l(\tau)=l(\sigma)+1$.  By multiplying the  equation $c_{i}\cdot  \sigma = \left(c_{i}\cdot  \tau\right)\cdot  (a,b)$ with $(c_i)^{-1}$ from the right we have $\sigma =\tau \cdot  (a,b)$.  	 This implies $\sigma\sma \tau$.
	
	"(ii) $\Leftrightarrow$ (iii)" Since $\eta \sma \mu $ if and only if $\eta^{-1} \sma \mu^{-1}$, and $\left(c_{j}\cdot  \sigma^{-1}\right)^{-1} = \sigma \cdot  d_{j} $, thus  $\left(c_{i}\cdot  \sigma^{-1}\right) \sma \left(c_{i}\cdot  \tau^{-1}\right)$ implies 
$ \left(\sigma \cdot  d_{i}\right)\sma \left(\tau \cdot  d_{i}\right)$ and vice versa. 
	\hfill $\Box$ 
\label{Bruhat_order_on_Blocks}
\end{subnum}

\begin{subnum}\normalfont{\textbf{Remark.}} Since $\mathbf{I}_{(i)}^{(n+1)} = c_i \cdot S_n$, every $w \in \mathbf{I}_{(i)}^{(n+1)}$ there exists $\tau \in S_n$ with $w = c_i \cdot \tau$. The proof of Lemma~\ref{Bruhat_order_on_Blocks}(2) shows that $(a,b)$ is a free-fall   of $\tau$ if and only if $(a,b)$ is a free-fall  of $w$, provided that $w\cdot (a,b) \in \mathbf{I}_{(i)}^{(n+1)}$. This  happens precisely when $a,b \neq n+1$. 
\begin{center}
$w = \begin{scriptsize} \left(
\begin{array}{ccccccc}
1 & \cdots  & a & \cdots & b &  \cdots & n+1\\
	\ast & \cdots & m & \cdots & r & \cdots & i 
\end{array}\right)
  \end{scriptsize}$, \   \  $w \cdot (a,b) = \begin{scriptsize} \left(
\begin{array}{ccccccc}
1 & \cdots  & a & \cdots & b &  \cdots & n+1\\
	\ast & \cdots & r & \cdots & m & \cdots & i 
\end{array}\right)
  \end{scriptsize} \in \textbf{I}_{(i)}^{(n+1)}$,
\end{center}
Furthermore, part (1) of Lemma~\ref{Bruhat_order_on_Blocks}  shows that  $c_{i+1} \cdot \tau $ is the unique small neighbor of $ c_i\cdot \tau =w$ within $\mathbf{I}_{(i+1)}^{(n+1)}$. It follows that $\left\{c_{i+1}\cdot \tau\right\} \cup \left\{c_i\cdot \sigma \mid \sigma \in \textbf{N}_{\textbf{s}}(\tau) \right\} \subseteq \textbf{N}_{\textbf{s}}(c_i\cdot \tau)$, and consequently $\left|\textbf{N}_{\textbf{s}}(c_i\cdot \tau)\right| \geq \left|\textbf{N}_{\textbf{s}}(c_i\cdot \tau)\right|+1$ as illustrated in the figure below.\\[2mm]

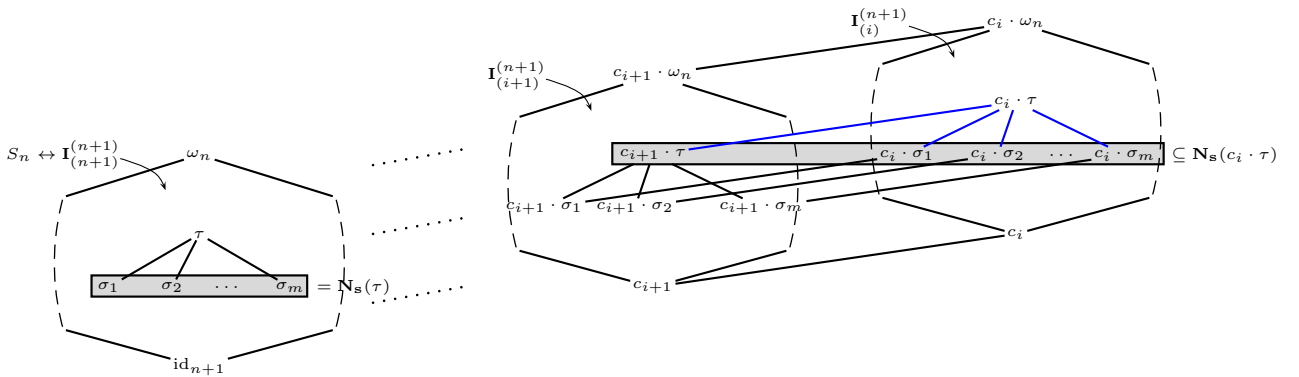
\begin{tiny}
\psset{xunit=1.2mm,yunit=0.7mm,runit=1mm}
\begin{pspicture}(0,0)(0,0)
\rput(17,-40){ \rput(-12,-12){\psframe*[linecolor=gray!30](0,0)(24,4.4)} \rput(-12,-12){\psframe[linecolor=black](0,0)(24,4.4)} \rput(17,-10){$=\textbf{N}_{\textbf{s}}(\tau)$}
\rput(19,0){\psline[linewidth=1pt,linestyle=dotted](0,0)(10,3)} 
\rput(19,13){\psline[linewidth=1pt,linestyle=dotted](0,0)(10,3)}
\rput(19,-13){\psline[linewidth=1pt,linestyle=dotted](0,0)(10,3)}  
\rput(-15,15){$S_n \leftrightarrow \textbf{I}_{(n+1)}^{(n+1)}$\pscurve[linewidth=0.4pt,linestyle=solid]{->}(0,0)(4,-3)(5,-6)} 
\rput(0,15){\rnode{omega}{$\omega_n$}}   \rput(20,20){\rnode{omegal}{}} 
\rput(-15,7){\rnode{g1}{}} \rput(15,7){\rnode{g2}{}}
\rput(0,-25){\rnode{id}{$\textrm{id}_{n+1}$}}
\rput(-15,-18){\rnode{m1}{}} \rput(15,-18){\rnode{m2}{}}
\rput(0,0){\rnode{tau}{$\tau$}}
\rput(-10,-10){\rnode{sigma_1}{$\sigma_1$}}
\rput(-3,-10){\rnode{sigma_2}{$\sigma_2$}} \rput(3,-10){$\ldots$}
\rput(10,-10){\rnode{sigma_m}{$\sigma_m$}}

\rput(50,15){ \rput(-4.6,-2){\psframe*[linecolor=gray!30](0,0)(61,4.4)} \rput(-4.6,-2){\psframe[linecolor=black](0,0)(61,4.4)}
\rput(-15,15){$\textbf{I}_{(i+1)}^{(n+1)}$\pscurve[linewidth=0.4pt,linestyle=solid]{->}(0,0)(4,-3)(5,-6)} 
\rput(0,15){\rnode{omegac}{$c_{i+1}\cdot \omega_n$}}    
\rput(-15,7){\rnode{g1c}{}} \rput(15,7){\rnode{g2c}{}}
\rput(0,-25){\rnode{idc}{$c_{i+1}$}}
\rput(-15,-18){\rnode{m1c}{}} \rput(15,-18){\rnode{m2c}{}}
\rput(0,0){\rnode{taus}{$c_{i+1}\cdot \tau$}}
\rput(-12,-10){\rnode{sigma_1s}{$c_{i+1}\cdot \sigma_1$}}
\rput(-2,-10){\rnode{sigma_2s}{$c_{i+1}\cdot \sigma_2$}} 
\rput(12,-10){\rnode{sigma_ms}{$c_{i+1}\cdot \sigma_m$}}
 }

\rput(90,25){  \rput(23,-10){$\subseteq\textbf{N}_{\textbf{s}}(c_{i}\cdot \tau)$}
\rput(-15,15){$\textbf{I}_{(i)}^{(n+1)}$\pscurve[linewidth=0.4pt,linestyle=solid]{->}(0,0)(4,-3)(5,-6)} 
\rput(0,15){\rnode{omegaci}{$c_{i}\cdot \omega_n$}}    
\rput(-15,7){\rnode{g1ci}{}} \rput(15,7){\rnode{g2ci}{}}
\rput(0,-25){\rnode{idci}{$c_{i}$}}
\rput(-15,-18){\rnode{m1ci}{}} \rput(15,-18){\rnode{m2ci}{}}
\rput(0,0){\rnode{tausi}{$c_{i}\cdot \tau$}} 
\rput(-12,-10){\rnode{sigma_1si}{$c_{i}\cdot \sigma_1$}}
\rput(-2,-10){\rnode{sigma_2si}{$c_{i}\cdot \sigma_2$}} \rput(5,-10){$\ldots$}
\rput(12,-10){\rnode{sigma_msi}{$c_{i}\cdot \sigma_m$}}
 }
 }

\psset{nodesep=1pt} 
\ncline{-}{omega}{g1} \ncline{-}{omega}{g2} \ncline{-}{id}{m1} \ncline{-}{id}{m2} \ncarc[linewidth=0.5pt,arcangle=15,linestyle=dashed]{-}{m1}{g1}
\ncarc[linewidth=0.5pt,arcangle=-15,linestyle=dashed]{-}{m2}{g2}

\ncline{-}{omegac}{g1c} \ncline{-}{omegac}{g2c} \ncline{-}{idc}{m1c} \ncline{-}{idc}{m2c} \ncarc[linewidth=0.5pt,arcangle=15,linestyle=dashed]{-}{m1c}{g1c}
\ncarc[linewidth=0.5pt,arcangle=-15,linestyle=dashed]{-}{m2c}{g2c}

\ncline{-}{omegaci}{g1ci} \ncline{-}{omegaci}{g2ci} \ncline{-}{idci}{m1ci} \ncline{-}{idci}{m2ci} \ncarc[linewidth=0.5pt,arcangle=15,linestyle=dashed]{-}{m1ci}{g1ci}
\ncarc[linewidth=0.5pt,arcangle=-15,linestyle=dashed]{-}{m2ci}{g2ci}

 \ncline{-}{tau}{sigma_1}\ncline{-}{tau}{sigma_2}\ncline{-}{tau}{sigma_m} 
 \ncline{-}{taus}{sigma_1s}\ncline{-}{taus}{sigma_2s}\ncline{-}{taus}{sigma_ms}
\ncline{-}{omegac}{omegaci}
 \ncline{-}{sigma_1si}{sigma_1s} \ncline{-}{sigma_2si}{sigma_2s} \ncline{-}{sigma_msi}{sigma_ms} \ncline{-}{idci}{idc}
\psset{nodesep=1pt, linecolor=blue}
\ncline{-}{taus}{tausi}
\ncline{-}{tausi}{sigma_1si}\ncline{-}{tausi}{sigma_2si}\ncline{-}{tausi}{sigma_msi}
\end{pspicture}
\end{tiny}

\textrm{  }\\[4cm]
\label{Number_of_neigh_on_Block}
\end{subnum}
This observation underlies the following

\begin{subnum}\normalfont{\textbf{Corollary.}} \textit{Let $\tau \in S_{n} = \textbf{I}_{(n+1)}^{(n+1)}$,  \  $c_i= s_{i}\cdot s_{i+1}\cdots s_{n}$ for $i\in \left[1,n\right]$ \  and \  $\Inv_{\textbf{f}}(c_i\cdot \tau)$ the set of free-fall inversions of $c_i\cdot \tau$.}

\begin{itemize}
\item[(1)]\textit{If  $\tau$ is critical  of degree $\textbf{\textit{i}}$, then $\tau$ is also critical permutation of degree $\textbf{\textit{i}}$  in $S_{n+1}$.}
	\item[(2)] \textit{$\left|\textbf{N}_{\textbf{s}}(c_i\cdot \tau)\right| =  \left|\textbf{N}_{\textbf{s}}( \tau)\right| +t_i$,  where $t_i= \left|\left\{(a,b)\in \Inv_{\textbf{f}}(c_i\cdot \tau) \mid  b=n+1\right\}\right| \geq 1$.}
	\item[(3)] \textit{If  \ $\left|\textbf{N}_{\textbf{s}}( \tau)\right|\geq n$, then $c_i\cdot  \tau$ is first degree critical permutation in $S_{n+1}$ for each $i\in \left[1,n+1\right]$.}
\end{itemize}

\textit{Proof.} (1)  Since $(a,n+1)\not\in \Inv(\tau)$, we have $\textbf{N}_{\textbf{s}}( \tau)\subseteq \textbf{I}_{(n+1)}^{(n+1)}$. By induction, we have  $\Lambda_{(\tau)} = \left\{\sigma \in S_{n+1}\mid \sigma\ma \tau\right\}\subseteq \textbf{I}_{(n+1)}^{(n+1)}$.   Thus  $\Lambda_{(\tau|_{S_n})} =\Lambda_{(\tau|_{S_{n+1}})}$. This shows that $\tau$, viewed as a permutation in  $S_{n}$ or in  $S_{n+1}$ has the same degree of "\textit{criticality}" (see ~\ref{isomorph}).

 (2) Obviously  $\Inv_{\textbf{f}}(c_i\cdot \tau) $ if the disjoint union of  $ \left\{(a,b)\in \Inv_{\textbf{f}}(c_i\cdot \tau) \mid a,b\neq n+1\right\} $ and  $ \left\{(a,b)\in \Inv_{\textbf{f}}(c_i\cdot \tau) \mid b= n+1\right\}$, while ~\ref{Number_of_neigh_on_Block} shows $ \left\{(a,b)\in \Inv_{\textbf{f}}(c_i\cdot \tau) \mid a,b\neq n+1\right\} =  \Inv_{\textbf{f}}(\tau)$. According to Lemma~\ref{neigh_lemma} we have $\left|\textbf{N}_{\textbf{s}}(c_i\cdot \tau)\right| = \left|\Inv_{\textbf{f}}(c_i\cdot \tau)\right| =  \left|\textbf{N}_{\textbf{s}}( \tau)\right| +t_i$ with $t_i :=\left|\left\{(a,b)\in \Inv_{\textbf{f}}(c_i\cdot \tau) \mid b= n+1\right\}\right|$. There is $d\in \left[1,n\right]$ with $(c_i\cdot \tau)(d) = i+1$, 
thus $(d,n+1)\in \Inv_{\textbf{f}}(c_i\cdot \tau) $. This implies $t_i\geq 1.$

(3) If $\left|\textbf{N}_{\textbf{s}}( \tau)\right|\geq n$, then $\tau = c_{n+1}\cdot  \tau $ is first degree critical  in $S_{n}$ and  in $S_{n+1}$  (see ~\ref{firs_degree_for_sure} and (1)).  By applying  (2),  we have that  $\left|\textbf{N}_{\textbf{s}}(c_{i}\cdot  \tau)\right| \geq n+1$ for all $i\in \left[1,n\right]$. Since $\left|\textbf{L}_{\textbf{t}}(c_{i}\cdot  \tau)\right| \leq \left|\left\{s \in S_{n+1} \mid l(s)=1\right\}\right| =n$ we have  $\left|\textbf{N}_{\textbf{s}}(c_{i}\cdot  \tau)\right| \neq \left|\textbf{L}_{\textbf{t}}(c_{i}\cdot  \tau)\right|$.
 \hfill $\Box$ 
	
	\label{if_larger_n_then_critical}
\end{subnum}

 As not every first-degree critical permutation $\tau\in S_n$ satisfies $\left|\textbf{N}_{\textbf{s}}(\tau)\right| \geq n$ a more detailed analysis is required.

\begin{numeris}\begin{center}$\textbf{N}_{\textbf{s}}(c_i\cdot\tau)$ \textbf{and} $\textbf{L}_{\textbf{t}}(c_i\cdot \tau)$\end{center}\end{numeris}
\begin{small}Recall that $\Red(\eta)$, the set of simple transpositions occurring in a reduced decomposition of $\eta$, coincides with the set of simple transpositions comparable with $\eta$, denoted by $\textbf{L}_{\textbf{t}}(\eta)$ (see Lemma~\ref{simple_transpositions_of_tau}(1)).
Furthermore, the set of free-fall inversions of $\eta$, denoted by $\Inv_{\textbf{f}}(\eta)$, has one-to-one correspondence  with the set of small neighbors of $\eta$ denoted by $\mathbf{N}_{\textbf{s}}(\eta)$: each free-fall inversion $(a,b)\in \Inv_{\textbf{f}}(\eta)$ gives rise to a unique small neighbor $\eta\cdot (a,b)\in \mathbf{N}_s(\eta)$ (see Lemma~\ref{neigh_lemma}).\end{small}\\

It is well known that whenever $\sigma \ma \tau$,  every reduced decomposition of $\sigma$ can be extended to a reduced decomposition of $\tau$, so that $\Red(\sigma) \subseteq \Red(\tau)$ (see \cite[Theorem 2.2.2]{Bjor}). 
Thus, for  $\sigma \ma \tau$ we always have $\left|\textbf{L}_{\textbf{t}}(\sigma)\right| \leq \left|\textbf{L}_{\textbf{t}}(\tau)\right|$. By contrast, the same monotonicity does not extend to the sets of  small neighbors. For example $(4231) \ma (4321)$, but $ \left|\mathbf{N}_{\textbf{s}}(4231) \right|= 4 > 3= \left|\mathbf{N}_{\textbf{s}}(4321)\right|$.

In the previous paragraph, we established that for each $\tau \in S_{n}=\textbf{I}_{(n+1)}^{(n+1)}$ and every $j\in \left[1,n\right]$ with $c_{j}=s_{j}\cdot s_{j+1}\cdots s_{n}$ we have  
$\tau = c_{n+1}\cdot \tau \sma c_{n}\cdot \tau  \sma c_{n-1}\cdot \tau \sma \cdots \sma c_{i+1}\cdot \tau \sma c_{i}\cdot \tau \sma \cdots \sma c_{1}\cdot \tau$
which implies $\left|\textbf{L}_{\textbf{t}}(c_{i}\cdot \tau)\right| \geq  \left|\textbf{L}_{\textbf{t}}(c_{i+1}\cdot \tau)\right|$.  However,   the analogous inequality for the sets of small neighbors, $\left|\mathbf{N}_{\textbf{s}}(c_{i}\cdot \tau )\right|\geq \left|\mathbf{N}_{\textbf{s}}(c_{i+1}\cdot \tau)\right|$ does not hold in general.

The principal statement of this paragraph, serving as a key  result for the proof of our main theorem, is the following:

\begin{subnum}\normalfont{\textbf{Lemma.}}  \textit{Let $\tau\in S_{n} \subseteq S_{n+1}$. For each $i\in \left[1,n+1\right]$  there is  $t_i\geq 0$  such that }
 \begin{center}
 \textit{$\left|\textbf{L}_{\textbf{t}}(c_{i}\cdot  \tau)\right| =\left|\textbf{L}_{\textbf{t}}( \tau)\right| +t_i$ \  \  \  \  and  \  \  \  \  
		 $\left|\textbf{N}_{\textbf{s}}(c_{i}\cdot  \tau)\right| \geq \left|\textbf{N}_{\textbf{s}}( \tau)\right| +t_i$.}
\end{center}
\label{kern_aussage}
\end{subnum}
Since the proof of this lemma relies on several preparatory observations 
concerning the structural differences between $\eta$ and $s_i \eta$, 
we postpone it to the end of the paragraph.

  Recall that 
every permutation 
 admits   a reduced decomposition that can be partitioned into \textit{blocks} via $\tau = \ger{b}_1 \cdot \ger{b}_2 \cdots \ger{b}_m$, where each \textit{block} $\ger{b}_i$ is a product of simple transpositions with consecutive indices, and for all $s_a\in \Red(\ger{b}_i)$ and $s_b\in \Red(\ger{b}_j)$ with $i<j$ we have $a<b$ (see Remark~\ref{Remark_about_blocks_general}). A permutation of this form can be represented schematically as follows:
\begin{center}
$\tau =  \begin{tiny}\Bigg(
\begin{array}{ccc}
	\cdots &  i & \cdots\\
	\cdots &  i & \cdots
\end{array}\underbrace{\framebox{$
\begin{array}{ccc}
	l_{1} & \cdots & k_{1} \\
	\ast & \cdots   & \ast 
\end{array}$
}}_{\ger{b}_1}
\begin{array}{ccc}
	\cdots &  i & \cdots\\
	\cdots &  i & \cdots
\end{array}\underbrace{\framebox{$
\begin{array}{ccc}
	l_{2} & \cdots & k_{2} \\
	\ast & \cdots   & \ast 
\end{array}$
}}_{\ger{b}_2} \cdots  \underbrace{\framebox{$
\begin{array}{ccc}
	l_{m} & \cdots & k_{m} \\
	\ast & \cdots   & \ast 
\end{array}$
}}_{\ger{b}_m}
\begin{array}{cccc}
	\cdots &  i & \cdots &n+1\\
	\cdots &  i & \cdots & n+1
\end{array}
\Bigg)\end{tiny} $
\end{center}

The proof of Lemma~\ref{kern_aussage} proceeds by induction on the number $m$ of blocks in the 
reduced decomposition of $\tau$. The determination of the integers $t_i$ for a single-block permutation $\tau = \ger{b}$  presents the core arguments, which extend inductively to multiple-block permutations.

 We consider $\tau = \ger{b} = s_{i_1}\cdot s_{i_2} \cdots s_{i_m}$ with $\Red(\ger{b}) =\left\{s_{i_1}, s_{i_2} ,\ldots, s_{i_m}\right\}  = \left\{s_{l}, s_{l+1} ,\ldots , s_{k-1}\right\}$. As previously observed,  $\tau(t)=t$ for all $t\in \left[1,n+1\right]$ with $t\not\in\left[l,k\right]$

 \begin{center}
 $\tau = \ger{b} = \begin{scriptsize}\Bigg(
 \begin{array}{ccc}
	1 & \cdots & l -1 \\
	1&\cdots&  l -1
\end{array}
\underbrace{\framebox{$
\begin{array}{cccc}
	l & l+1 & \cdots &k \\
	\ast & \ast & \cdots &\ast
\end{array}$
}}_{\ger{b}}
\begin{array}{ccc}
	k+1 & \cdots & n\\
	k+1 & \cdots & n
\end{array}
\begin{array}{c}
	n+1\\
	n+1
\end{array} \Bigg)\end{scriptsize}$. 
\end{center}
 Note  $\left[1,n\right] =\left[1,l-1\right] \dot{\cup}  \left[l,k-1\right] \dot{\cup}  \left[k,n\right]$, and  for $c_j\cdot \tau = s_{j}\cdot s_{j+1}\cdots s_{n} \cdot s_{i_1}\cdot s_{i_2} \cdots s_{i_m}$ we have
\begin{center}
$\Red(c_j\cdot \tau) = \Red(\tau) \cup \Red(c_j) = \left\{s_{l}, s_{l+1} ,\ldots , s_{k-1}, s_{j},s_{j+1}, \ldots , s_{n-1}, s_{n}\right\}$.
\end{center}
 Since $\Red(\tau)\cap \Red(c_i) =\emptyset$ for $i\in\left[k,n\right] $,  $\Red(\tau)\cap \Red(c_i) =\left\{s_{i}, s_{i+1} ,\ldots , s_{k-1}\right\}$ for $i\in\left[l,k-1\right] $,  and  $\Red(\tau)\cap \Red(c_i) =\Red(\tau)$ for $i\in\left[1,l-1\right]$, we have
\begin{center}
$
\begin{array}{cclcl}
	\Red(c_i\cdot \tau)& =& \Red(\tau) \dot{\cup} \left\{s_{i},s_{i+1}, \ldots , s_{n-1}, s_{n}\right\}  & \textrm{if} & i\in\left[k,n\right],\\[1mm]
	\Red(c_i\cdot \tau)& =& \Red(\tau) \dot{\cup} \left\{ s_{k},s_{k+1}, \ldots , s_{n-1}, s_{n}\right\}  & \textrm{if} & i\in\left[l,k-1\right],\\[1mm]
	\Red(c_i\cdot \tau)& =& \Red(\tau) \dot{\cup} \left\{s_{i},s_{i+1}, \ldots ,   s_{l-1}, s_{k},\ldots , s_{n-1}, s_{n}\right\}  & \textrm{if} & i\in\left[1,l-1\right].
\end{array}
$
\end{center}
By Lemma~\ref{simple_transpositions_of_tau}(1),   we have   $\left|\textbf{L}_{\textbf{t}}(c_i\cdot \tau)\right| =  \left|\textbf{L}_{\textbf{t}}( \tau)\right| +t_i$,  for  $t_i= \begin{scriptsize}\left\{
\begin{array}{ccl}
	n-i+1& if & i\in \left[k,n\right],\\
	n-k+1 &if& i\in \left[l,k-1\right],\\
	n-k +1 +l-i &if & i\in \left[1,l-1\right].
\end{array}
\right.\end{scriptsize}$ We now proceed to show that $\left|\textbf{N}_{\textbf{s}}(c_i\cdot \tau)\right| \geq   \left|\textbf{N}_{\textbf{s}}( \tau)\right| +t_i$.
\begin{itemize}
	\item If $i\in \left[k,n\right]$ we have  $s_i\not\in \Red(c_{i+1}\cdot \tau)$. Since $c_{i}\cdot  \tau = s_i\cdot c_{i+1}\cdot  \tau$,  Lemma~\ref{simple_transpositions_of_tau}(2) yields 
 $\left|\textbf{N}_{\textbf{s}}(c_{i}\cdot  \tau)\right| = \left|\textbf{N}_{\textbf{s}}(c_{i+1}\cdot  \tau)\right| +1$. 
Applying this relation inductively, we obtain  $\left|\textbf{N}_{\textbf{s}}(c_{i}\cdot  \tau)\right| =\left|\textbf{N}_{\textbf{s}}( c_{i+1}\cdot\tau)\right| +1 = \left|\textbf{N}_{\textbf{s}}( c_{i+2}\cdot\tau)\right| +2 = \cdots = \left|\textbf{N}_{\textbf{s}}( c_{i +(n+1-i)}\cdot\tau)\right| +n+1-i = \left|\textbf{N}_{\textbf{s}}( \tau)\right|  +t_i$.

	\item Let $i\in \left[l,k-1\right]$.  Since right-multiplication of a permutation by $s_j$ swaps the entries $j$ and $j+1$ in the bottom line, applying the multiplication of $s_n, s_{n-1}, \dots, s_i$ iteratively on $\tau $ we obtain that the entries in the bottom line of $c_i \cdot \tau = s_i\cdot s_{i+1}  \cdots s_{k-1}\cdot s_{k}\cdot s_{k+1} \cdots s_{n}\cdot \tau $  are arranged according to the following distribution
\begin{center}
$c_{i}\cdot\tau\!=\! \begin{scriptsize}\Bigg(
\begin{array}{ccc}
	1 & \cdots & l -1 \\
	1&\cdots&  l -1
\end{array}
\framebox{$
\begin{array}{ccccc}
	l& \cdots &m_i&\cdots &k  \\
	\ast & \cdots & i+1&\cdots& \ast 
\end{array}$
}
\begin{array}{ccccc}
	k+1 &  \cdots &j&   \cdots & n\\
	k+2 &  \cdots &j+1&  \cdots & n+1
\end{array}\begin{array}{c}
	n+1\\
	 i
\end{array}\Bigg)\end{scriptsize}$.
\end{center}
Obviously,  $\left\{(m_i, n+1)\right\}\cup \left\{ (a,n+1) \mid  a\in \left[k+1,n\right]\right\} 
\subseteq \left\{(a,b)\in \Inv_{\textbf{f}}(c_i\cdot\tau) \mid b=n+1\right\}$. Thus $\left|\left\{(a,b)\in \Inv_{\textbf{f}}(c_i\cdot\tau) \mid b=n+1\right\}\right| \geq \left|\left\{(m_i, n+1)\right\}\cup \left\{ (a,n+1) \mid  a\in \left[k+1,n\right]\right\} \right| = n-k+1$. Based on Corollary~\ref{if_larger_n_then_critical} we have 
$\left|\textbf{N}_{\textbf{s}}(c_i\cdot \tau)\right| \geq   \left|\textbf{N}_{\textbf{s}}( \tau)\right| +n-k+1 = \left|\textbf{N}_{\textbf{s}}( \tau)\right| +t_i$.

	\item For  $i\in \left[1,l-1\right]$,  we have  $s_i\not\in \Red(c_{i+1}\cdot \tau)$. Since $c_{i}\cdot  \tau = s_i\cdot c_{i+1}\cdot  \tau$,  we have 
 $\left|\textbf{N}_{\textbf{s}}(c_{i}\cdot  \tau)\right| = \left|\textbf{N}_{\textbf{s}}(c_{i+1}\cdot  \tau)\right| +1$ (see ~\ref{simple_transpositions_of_tau}(2)). Inductively, we obtain  $\left|\textbf{N}_{\textbf{s}}(c_{i}\cdot  \tau)\right| =\left|\textbf{N}_{\textbf{s}}( c_{i+1}\cdot\tau)\right| +1 = \left|\textbf{N}_{\textbf{s}}( c_{i+2}\cdot\tau)\right| +2 = \cdots = \left|\textbf{N}_{\textbf{s}}( c_{l}\cdot\tau)\right| +l-i$. Since  $\left|\textbf{N}_{\textbf{s}}( c_{l}\cdot\tau)\right|\geq \left|\textbf{N}_{\textbf{s}}(\tau)\right|+n-k+1$, we obtain $\left|\textbf{N}_{\textbf{s}}(c_{i}\cdot  \tau)\right| \geq \left|\textbf{N}_{\textbf{s}}(\tau)\right|+n-k+1 +l-i$. Thus $\left|\textbf{N}_{\textbf{s}}(c_{i}\cdot  \tau)\right| \geq \left|\textbf{N}_{\textbf{s}}(\tau)\right|+t_i$.
\end{itemize}

	\textit{Proof of Lemma~\ref{kern_aussage}.} Consider a reduced decomposition of  $\tau = \ger{b}_1\cdot \ger{b}_2\cdots \ger{b}_m \in S_n=\textbf{I}_{(n+1)}^{(n+1)}$ as described above, and for each block $\ger{b}_j$ we have  $\Red(\ger{b}_j)=\left\{s_{i} \mid i\in \left[l_{j},   k_{j}-1\right]\right\} $ with $1\leq l_1<k_1<l_2<k_2 <\cdots <l_m<k_m\leq n$. By induction on $m$ we have  $\left|\textbf{L}_{\textbf{t}}(c_{i}\cdot  \tau)\right| =\left|\textbf{L}_{\textbf{t}}( \tau)\right| +t_i$ \  and 
		 $\left|\textbf{N}_{\textbf{s}}(c_{i}\cdot  \tau)\right| \geq \left|\textbf{N}_{\textbf{s}}( \tau)\right| +t_i$ for 
	\begin{center}
	$t_i = \begin{small}\left\{
		\begin{array}{ccl}
	n-i+1 & \textrm{if} & i\in \left[k_m, n\right] \\
		t_{k_j} & \textrm{if} & i\in \left[l_{j}, k_{j}-1\right] \ \textrm{and} \ 1\leq j\leq m \\
		t_{k_j} +l_{j}- i & \textrm{if} & i\in \left[k_{j-1}, l_{j}-1\right]\ \textrm{and} \ 1\leq j\leq m \\
		t_{k_1} +l_{1}- i & \textrm{if} & i\in \left[1, l_{1}-1\right]
	\end{array}\right.\end{small}
	 $
\end{center}
For $m=1$,  the claim holds, as shown previously (note that $t_i = t_{k_1} = n-k_1+1$ for all $i \in [l_1,k_{1}-1]$ and $t_i = t_{k_1} +l_{1}- i  = n-k_1+1  +l_1-i$ for all $i\in \left[1, l_{1}-1\right]$.)  For the induction step, we consider the intervals $\left[k_{j-1}, l_{j}-1\right]$, $\left[l_{j}, k_{j}-1\right]$, $\left[k_{j}, l_{j+1}-1\right]$ and apply the same structural properties of permutations with this type of reduced decomposition.
\hfill $\Box$

\section{\textbf{Proof of Theorem}~\ref{Main_Theorem}}
\label{proof_main_Theorem}

Combining Lemmas~\ref{simple_transpositions_of_tau} and ~\ref{kern_aussage}, we obtain the following 

	\begin{numeris}\normalfont{\textbf{Corollary}} \textit{If $\tau$ is first degree critical permutation in $S_n$, then $c_i\cdot  \tau$ 
is  first degree critical permutation in $S_{n+1}$ for every  $i\in \left[1, n+1\right]$.}\\ 

\textit{Proof.} If $\tau$ is first-degree critical, then  $\left|\textbf{N}_{\textbf{s}}( \tau)\right| > \left|\textbf{L}_{\textbf{t}}( \tau)\right|$ (see ~\ref{simple_transpositions_of_tau}(4)). Based on Lemma~\ref{kern_aussage} we have 
$\left|\textbf{N}_{\textbf{s}}(c_i\cdot  \tau)\right| \geq \left|\textbf{N}_{\textbf{s}}( \tau)\right| + t_i > \left|\textbf{L}_{\textbf{t}}( \tau)\right| + t_i = \left|\textbf{L}_{\textbf{t}}( c_i\cdot \tau)\right| $. Hence, for every $i \in [1,n+1]$, the permutation $c_i \cdot \tau$ 
is also first-degree critical in $S_{n+1}$   (by ~\ref{simple_transpositions_of_tau}(4)). 
\hfill $\Box$ 
	\label{proof_of_the_second_part}
\end{numeris}

 \begin{numeris}\normalfont{\textbf{Remark.}} If the Bruhat-diagrams of $\left(\Lambda_{(\tau)}, \Bruh \right)$ and $\left(\Lambda_{(\widetilde{\tau})}, \Bruh \right)$  have the same structural characteristic, then  "\textit{criticality status}" of $\tau$ and $\widetilde{\tau}$ is interconnected, because it implies that for all $\textbf{k}$ we have  $\left|\left\{\eta\in \Lambda_{(\tau)} \mid l(\eta)=\textbf{k}\right\}\right| = \left|\left\{\eta\in \Lambda_{(\widetilde{\tau})} \mid l(\eta)=\textbf{k}\right\}\right| $. In other words: if there exists a bijective, order-preserving function $f:\Lambda_{(\tau)} \to \Lambda_{(\widetilde{\tau})}$  (i.e. $\sigma \Bruh \mu$ if and only if $f(\sigma)\Bruh f(\mu)$ for all $\sigma,\mu\in \Lambda_{(\tau)}$), then we say $\left(\Lambda_{(\tau)}, \Bruh \right)$ and $\left(\Lambda_{(\widetilde{\tau})}, \Bruh \right)$ \textit{isomorphic} and use the notation $\Lambda_{(\tau)} \cong \Lambda_{(\widetilde{\tau})}$. Note, that in the Bruhat-diagram  $\sigma $ and $f(\sigma)$ are placed at the same level on the length-scale,   and   $\sigma $ is connected to $ \mu$ by a line if and only if $f(\sigma) $ is connected to $f(\mu)$ as well. Consequently, $\Lambda_{(\tau)} \cong \Lambda_{(\widetilde{\tau})}$ implies that $\tau$ is a critical permutation of degree $\textbf{\textit{i}}$ if and only if $\widetilde{\tau}$ is. 
 \label{isomorph}
\end{numeris}

\begin{numeris}\normalfont{\textbf{Lemma.}} 	\textit{Let $\sigma \in S_m$ and $\omega_m=(m\ldots 321) $ the  longest element. The
following statements are equivalent:}
\begin{itemize}
	\item[(i)] \textit{$\sigma$ is a critical permutation of degree $\textbf{\textit{i}}$.}
	\item[(ii)] \textit{$\omega_m \cdot  \sigma \cdot  \omega_m$ is a critical permutation of degree $\textbf{\textit{i}}$.}
	\item[(iii)] \textit{$\sigma^{-1}$ is a critical permutation of degree $\textbf{\textit{i}}$.}
\end{itemize}
\textit{Proof.} The mappings $S_m \to S_m$ via $\textbf{x} \mapsto \omega_m\cdot  \textbf{x} \cdot  \omega_{m}$ as well as via $\textbf{x} \mapsto \textbf{x}^{-1}$ are inner group automorphisms and  Bruhat order automorphisms   in the sense of the following being equivalent 
  \begin{center}
	(i) $\sigma \sma \tau $, \  \      (ii)  $\sigma^{-1} \sma \tau^{-1} $  \  \    (iii) $\left(\omega_n\cdot  \sigma \cdot  \omega_m \right) \sma \left(\omega_m\cdot \tau \cdot  \omega_m\right) $   \  \    (iv) $\left(\omega_m\cdot  \sigma^{-1} \cdot  \omega_m \right) \sma \left(\omega_m\cdot \tau^{-1} \cdot  \omega_m\right)$ 
\end{center}
	 (see  \cite[Proposition 3.1.5 ]{Bjor}). 
 Thus, the functions $f:\Lambda_{(\tau)} \to \Lambda_{(\omega_m\cdot \tau\cdot \omega_m)}$ with $f(\sigma)= \omega_m\cdot \sigma\cdot \omega_m$,  and $g:\Lambda_{(\tau)} \to \Lambda_{(\tau^{-1})}$ with $g(\sigma) = \sigma^{-1}$ are bijective and order-preserving functions, so 
\begin{center}
$\Lambda_{(\tau)}  \cong \Lambda_{(\omega_m\cdot \tau\cdot \omega_m)} \cong \Lambda_{(\tau^{-1})} \cong \Lambda_{(\omega_m\cdot \tau^{-1}\cdot \omega_m)}$ \  \  for each $\tau\in S_m$
\end{center}
implies the prove of the equivalent statements.      \hfill $\Box$ 
\label{critical_if_and_only_if}
\end{numeris}

This lemma constitutes the final observation, completing the proof of the  Theorem~\ref{Main_Theorem}, which asserts that if  $\tau\in S_n$ is first degree critical permutation, then for every $i\in \left[1, n+1\right]$, permutations
$c_i\cdot \tau$,   $c_i\cdot \tau^{-1}$,   $\tau\cdot d_i$,  $\tau^{-1}\cdot d_i$,  $\omega_{n+1}\cdot c_i\cdot \tau \cdot \omega_{n+1}$, $\omega_{n+1}\cdot c_i\cdot \tau^{-1} \cdot \omega_{n+1}$,   $\omega_{n+1}\cdot d_i \cdot \tau \cdot \omega_{n+1}$, and  $\omega_{n+1}\cdot d_i \cdot \tau^{-1} \cdot \omega_{n+1}$ are first degree critical permutations in $S_{n+1}$.\\

\textit{Proof of the Theorem~\ref{Main_Theorem}.} A permutation $\tau\in S_n$ is first-degree critical if and only if $\tau^{-1}\in S_n$
is. Based on Corollary~\ref{proof_of_the_second_part}  for every $i\in \left[1, n+1\right]$, permutations 
$c_i\cdot \tau$ and   $c_i\cdot \tau^{-1}$ are critical in first degree in $S_{n+1}$. Based on ~\ref{critical_if_and_only_if} we have that $\left(c_i\cdot \tau\right)^{-1} = \tau^{-1}\cdot d_i$ and   $\left(c_i\cdot \tau^{-1}\right)^{-1}=\tau\cdot d_i$, as well as their $\omega_{n+1}$-conjugates are first degree critical permutations as well. \hfill $\Box$

\begin{subnum}\normalfont{\textbf{Remark.}}	According to Jantzen's computations, among the $32$ permutations in $S_5$ yielding non-simple Verma multiplicities, there are \textit{four} elements, $w$,  with $\left[M(\id_5):L(w)\right]=3$ and \textit{one}  with $\left[M(\id_5):L(w)\right]=4$ (see \cite[5.24]{Jantzen}). Our analysis of the critical permutations in $S_5$ shows that there are  \textit{five} permutations $\tau$ for which  
$\textbf{N}_{\textbf{s}}(\tau)- \textbf{L}_{\textbf{t}}(\tau) = 2$, one second-degree critical permutation, and the remaining  permutations satisfy  $\textbf{N}_{\textbf{s}}(\tau)- \textbf{L}_{\textbf{t}}(\tau) = 1$. These observations indicate that there is no evident relationship between   $\textbf{N}_{\textbf{s}}(\tau)- \textbf{L}_{\textbf{t}}(\tau)$ and  Verma multiplicities.
\label{what_we_absurved_until_now}
\end{subnum}

\section{Other Classes of Critical Permutations}
\begin{small}As noted earlier, not every first-degree critical permutation in $S_{n+1}$ arises from 
a first-degree critical permutation in $S_n$. 
 In this section,  we describe several families of permutations with special forms 
that in general do not arise 
 from the construction described in the Theorem~\ref{Main_Theorem}. This perspective allows us to investigate additional structures and patterns, 
yielding a more complete understanding of critical permutations.  Since the arguments rely primarily on counting inversions and are straightforward, 
we leave out the full proofs, providing instead an outline of the main strategies and key ideas.\end{small}

\begin{numeris}\normalfont{\textbf{First-degree critical permutations}}.  By Lemma~\ref{if_larger_n_then_critical}(3), permutations $\tau\in S_n$ satisfying 
$\left|\mathbf{N}_{\textbf{s}}(\tau)\right| \geq n$ 
are first-degree critical. The following permutations fall into this class.

\begin{itemize}
	\item[1.]  Transposition  $\tau = (i,i+k)$, with $k\geq \frac{n +2}{2}$ is first-degree critical permutation.
	
  For every $k$ we have $\Inv_{\textbf{f}}(i,i+k) =  \Inv(i,i+k) \backslash \left\{(i,i+k)\right\}$, and $\left| \Inv(i,i+k)\right| = 2\cdot k -1$, thus 	$\left|\textbf{N}_{\textbf{s}}(i,i+k)\right| \stackrel{~\ref{neigh_lemma}}{=} \left|\Inv_{\textbf{f}}(i,i+k)\right| =  2\cdot k -2 \geq n $ for $k\geq \frac{n +2}{2}$.

	\item[2.] If  $n\geq 4$, then  $\tau = \begin{scriptsize}\left(
\begin{array}{cccccccccccc}
	1& 2& \cdots & l & \cdots & n-k & n-k+1 & n-k+2 & \cdots & n-k+l & \cdots & n \\
	k+1 & k+2 & \cdots & k+l &  \cdots & n & 1 & 2 & \cdots &l & \cdots & k
\end{array}
\right)\end{scriptsize}$  is first-degree critical permutation  for every   $ k \in \left[2, n-2\right]$.

It can be established that $\Inv_{\textbf{f}}(\tau)= \Inv(\tau) = \left\{(t,s) \mid t\in \left[1,n-k\right], \  s\in \left[n-k+1,n\right] \right\}$. Thus 
 $\left|\textbf{N}_{\textbf{s}}(\tau)\right| = \left|\Inv(\tau)\right| =  k \cdot (n-k) $. For all $a,b\in \NN$ with $a,b\geq 2$ we have $a\cdot b \geq a+b$. Since $2\leq k\leq n-2$ we have $k,n-k\geq 2$, thus $\left|\textbf{N}_{\textbf{s}}(\tau)\right| = k\cdot (n-k) \geq k + n-k = n$.

	\item[3.]  If  $a\geq 3$, then 
	$\tau = \begin{scriptsize}\left(
\begin{array}{ccccccccccccccc}
	1 &2& \cdots & a-k & \cdots &a -2 &  a -1  & a & a+1 & a+2 & \cdots & a+ k & \cdots & n \\
	n-1 & n-3 & \cdots & 2 k & \cdots &4 &  2 & n & 1 & 3 & \cdots & 2 k-1 & \cdots & n-2
\end{array}
\right)\end{scriptsize}$ is first-degree critical  for all $n = 2a -1$

	Analogously to the previous cases, a complete determination of the set $\Inv_{\textbf{f}}(\tau)$ shows that 	
	 $\left|\textbf{N}_{\textbf{s}}(\tau)\right| = \frac{3\cdot n -5}{2}\geq n$ for all $n\geq 5$.  Since $\left|\Inv_{\textbf{f}}(\tau)\right| = \left|\Inv_{\textbf{f}}(\omega_{n} \cdot  \eta \cdot  \omega_{n})\right|$ we have that $\tau$ and 
	 $\omega_{n} \cdot  \eta \cdot  \omega_{n} = \left(\framebox{ 3  \ 5 \ldots   n-2 \ n} \: 1\: \framebox{n-1 \  n-3   \ldots \ 4 \ 2}\right)$ are critical in $S_n$. 	
	\item[4.] Let  $\tau = (1,n)\cdot  s_{i_1} \cdot  s_{i_2} \cdots  s_{i_{k}} $ where $2\leq i_j\leq n-2$ and $i_{j}\neq i_{j'}$ for $j\neq j'$. Then $\left|\textbf{N}_{\textbf{s}}(\tau)\right| = k+4$. Hence, for $n \geq 5$ and $k \geq n-4$, 
permutations of this form are first-degree critical.
\end{itemize}
\label{citical_2}
\end{numeris}
\begin{numeris}\normalfont{\textbf{"\textit{Second-degree}" critical permutation}}.  For  $n= 2a-1$ with $a\geq 3$, the permutation 
\begin{center}
$p_n :=\begin{scriptsize}  \left(
\begin{array}{ccccccccc}
	1 & 2 & \cdots & a-1 & a & a+1& a+2 & \cdots & n\\
	a+1 & a+2 & \cdots & n & a & 1& 2 & \cdots & a-1
\end{array}
\right) \end{scriptsize}$
\end{center}
is self-inverse and its own  $\omega_n$-conjugate. Thus for every $\sigma\in S_n$ with $\sigma\ma p_n$ we have $\sigma^{-1}\ma p_n$ and $\omega_n\cdot \sigma \cdot \omega_n \ma p_n$. Using these properties, it can be verified that that $p_n$ is comparable with all $n-1$ simple transpositions in $S_{n}$, and  we can identify all transpositions $(i,j)$ of $p_n$ 
that yield \textit{small neighbors} $(i,j)\cdot p_n$,
as illustrated in the figure below.\\

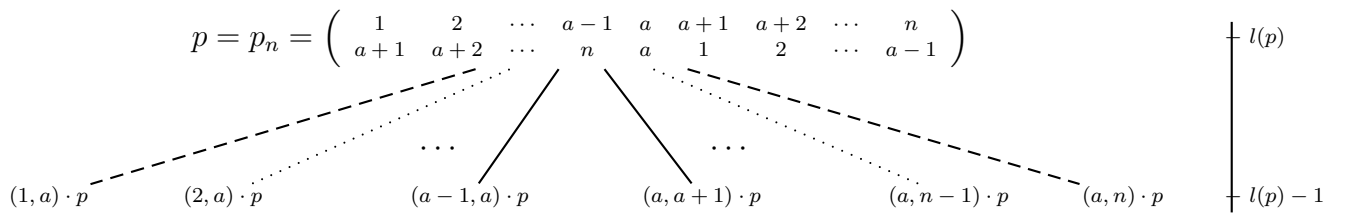
\psset{xunit=1.92mm,yunit=2.1mm,runit=1mm}
\begin{pspicture}(0,0)(0,0)
\rput(38,0){  
\rput(0,0){\rnode{top}{$p= p_n= \begin{scriptsize}  \left(
\begin{array}{ccccccccc}
	1 & 2 & \cdots & a-1 & a & a+1& a+2 & \cdots & n\\
	a+1 & a+2 & \cdots & n & a & 1& 2 & \cdots & a-1
\end{array}
\right) \end{scriptsize}$}
  }
\rput(-5.5,0){\rput(50,1){\rnode{sk0}{}}	 \rput(50,-11){\rnode{sku}{}}   \rput[l](49.5,0){\begin{scriptsize}$ - \      l(p)$\end{scriptsize}}	
 \rput[l](49.5,-10){\begin{scriptsize}$- \   l(p) -1 $\end{scriptsize}}} 
\rput(37,-10){\rnode{a1}{\begin{scriptsize}$(a,n) \cdot  p $\end{scriptsize}}}	
	\rput(25,-10){\rnode{a2}{\begin{scriptsize}$(a,n-1) \cdot  p $\end{scriptsize}}}		
	\rput(8,-10){\rnode{an-1}{\begin{scriptsize}$(a,a+1) \cdot  p $\end{scriptsize}}}		
	\rput(-37,-10){\rnode{1a}{\begin{scriptsize}$(1,a) \cdot  p $\end{scriptsize}}}		
	\rput(-25,-10){\rnode{2a}{\begin{scriptsize}$(2,a) \cdot  p $\end{scriptsize}}}	
	\rput(-8,-10){\rnode{a-1a}{\begin{scriptsize}$(a-1,a) \cdot  p $\end{scriptsize}}}	
	\rput(-10,-7){\rnode{pointl}{$\cdots$}}	
	\rput(10,-7){\rnode{pointr}{$\cdots$}}	
	}
	\ncline{sk0}{sku}
	\psset{nodesep=1pt,arrows=-}
	\ncline[linestyle=dashed]{top}{a1}
	\ncline[linestyle=dotted]{top}{a2}
	\ncline{top}{an-1}
	\ncline[linestyle=dashed]{top}{1a}
	\ncline[linestyle=dotted]{top}{2a}
	\ncline{top}{a-1a}
\end{pspicture}

\textrm{  }\\[2cm]	

Using mathematical induction, it is easy to show that the number of permutations $\sigma$ in $S_{n}$ with $l(\sigma) = 2$ is  $\left|\left\{\sigma\in S_n \mid l(\sigma)= 2\right\}\right| = \frac{(n+1)(n-2)}{2} $. By analyzing the inversions of the small neighbors of $p$ and using $p=p^{-1}=\omega_{n}\cdot p\cdot \omega_{n}$ we can determine that the number of permutations in $\left\{\sigma\in \Lambda_{(p)} \mid l(\sigma)= l(p)-\textbf{2}\right\}$ exceeds $\frac{(n+1)(n-2)}{2}$. 
These arguments form the basis for establishing 
\begin{itemize}
	\item $\left|\left\{\sigma\in \Lambda_{(p)} \mid l(\sigma)=l(p)- 1\right\}\right| = \left|\left\{\sigma\in \Lambda_{(p)} \mid l(\sigma)=1\right\}\right|$.
	\item $\left|\left\{\sigma\in \Lambda_{(p)} \mid l(\sigma)=l(p)- 2\right\}\right| \neq \left|\left\{\sigma\in \Lambda_{(p)} \mid l(\sigma)=2\right\}\right|$.
	\end{itemize}
and consequently, demonstrate that $p$ is a  critical permutation of degree \textit{\textbf{2}}.
\label{second_degree}
\end{numeris}
\label{sec:OtherCriticalPermutations}

\textbf{Acknowledgments.}  
I would like to thank Ahmed Elbosari for  computing and classifying all critical permutations in $S_n$ for $n\leq 7$.  This work clarified  how critical permutation of $S_n$ give rise to critical permutations in  $S_{n+1}$  and facilitated the identification of their relevant structures.

\begin{small}

\end{small}

\end{document}